# A SMALL-GAIN THEOREM FOR A WIDE CLASS OF FEEDBACK SYSTEMS WITH CONTROL APPLICATIONS

**Iasson Karafyllis[*] and Zhong-Ping Jiang[1,**]**

[*]Department of Environmental Engineering, Technical University of Crete, 73100, Chania, Greece
email: ikarafyl@enveng.tuc.gr

[**]Department of Electrical and Computer Engineering, Polytechnic University, Six Metrotech Center, Brooklyn, NY 11201, U.S.A.
email: zjiang@control.poly.edu

**Abstract**
A Small-Gain Theorem, which can be applied to a wide class of systems that includes systems satisfying the weak semigroup property, is presented in the present work. The result generalizes all existing results in the literature and exploits notions of weighted, uniform and non-uniform Input-to-Output Stability (IOS) property. Applications to partial state feedback stabilization problems with sampled-data feedback applied with zero order hold and positive sampling rate, are also presented.

**Keywords:** Input-to-Output Stability, Control Systems, Small-Gain Theorem.

## 1. Introduction

A common feature of stability analysis for complex interconnected systems is the application of small-gain results. Small-gain theorems for continuous-time finite-dimensional systems expressed in terms of "nonlinear gain functions" have a long history (see [10,27,44,45] and the references therein). A nonlinear small-gain result was presented in [10], which allowed numerous applications to feedback stabilization problems. The methodology presented in [10] was followed by many researchers (see [11,12,18,20,40,43]). A common characteristic of current research on nonlinear small-gain results in Mathematical Systems Theory is the use of the notion of uniform Input-to-State Stability (ISS), introduced by E.D. Sontag in [37] for systems described by Ordinary Differential Equations, or the notion of uniform Input-to-Output Stability (IOS), introduced by E.D. Sontag and Y. Wang in [39] (also see [10]) and extended in [8]. Small-Gain theorems for discrete-time systems can be found in [13,14,15].

Extensions of small-gain results were presented recently in the literature. A non-uniform in time small-gain theorem for continuous-time finite-dimensional systems was presented in [18]. Moreover, in [20,40] small-gain results for wide classes of systems were provided. The classes of systems considered in [20,40] satisfy the classical "semigroup property" (see [20,21,22,38]). Small-gain results for hybrid systems satisfying the classical "semigroup property" were recently presented in [26].

An important feature of certain hybrid systems is that they do not satisfy the classical "semigroup property": for example, the solution $x(t)$ of a system $\Sigma$ with initial condition $x(t_0) = x_0$ does not coincide (in general) for $t \geq t_1 > t_0$ with the solution $\tilde{x}(t)$ of $\Sigma$ with initial condition $\tilde{x}(t_1) = x(t_1)$. Such systems arise when sampled-data feedback laws are applied to finite-dimensional control systems or when numerical discretization schemes are applied for the numerical solution of a system of ordinary differential equations. However, from a mathematical point of view, these structures cannot be considered as "systems" in the sense given in [16,20,38]. This feature has important consequences, since the researcher cannot use the tools developed for systems theory and mathematical control theory. In [21,22] the notion of a system was relaxed so that the "semigroup property" does not hold in a strict sense.

---
[1] Corresponding author

Moreover, the modification introduced allows the results obtained in [20] to hold. Thus we are in a position to develop a complete stability theory, which covers the systems that satisfy the so-called "weak semigroup property".

The purpose of the present work is to present a small-gain result (Theorem 3.1 and Corollary 3.4), which

∗ can be applied to a very general class of systems (including systems that do not satisfy the classical "semigroup" property).
∗ unifies all existing results, which make use of uniform or non-uniform and weighted notions of Input-to-State Stability (ISS) or Input-to-Output Stability (IOS).
∗ can be used directly for the solution of sampled-data feedback stabilization problems or problems of numerical stability of discretization schemes.
∗ can be applied to uncertain time-varying systems with vanishing or non-vanishing perturbations.

We believe that the main result of the present work is a valuable tool for establishing stability and will be used frequently in future research. However, we would like to emphasize the theoretical significance of our main result: it is a method for establishing qualitative properties expressed in a very general framework that unifies works from various fields as well as different stability notions. The results presented in the paper can be extended without much difficulty to the case of local stability notions.

The contents of this paper are presented as follows. In Section 2 we provide the definitions of the notions used and several examples of systems that have the "Boundedness-Implies-Continuation" (BIC) property. In Section 3 the main result is stated and proved. In Section 4, it is shown how the main result of the present work can be applied to an ISS stabilization problem of a certain class of systems with partial-state sampled-data feedback. It should be emphasized that sampled-data control systems cannot be handled with Small-Gain results that have appeared so far in the literature, since sampled-data control systems do not satisfy the classical semigroup property. Finally, Section 5 contains the conclusions of the paper. The proofs of some basic results are given in the Appendix.

**Notations** Throughout this paper we adopt the following notations:
∗ We denote by $K^+$ the class of positive, continuous functions defined on $\Re^+$. We say that a function $\rho : \Re^+ \to \Re^+$ is of class $N$, if $\rho$ is continuous, non-decreasing with $\rho(0) = 0$. By $K$ we denote the set of positive definite, increasing and continuous functions. We say that a positive definite, increasing and continuous function $\rho : \Re^+ \to \Re^+$ is of class $K_\infty$ if $\lim_{s \to +\infty} \rho(s) = +\infty$. By $KL$ we denote the set of all continuous functions $\sigma = \sigma(s,t) : \Re^+ \times \Re^+ \to \Re^+$ with the properties: (i) for each $t \geq 0$ the mapping $\sigma(\cdot, t)$ is of class $K$; (ii) for each $s \geq 0$, the mapping $\sigma(s, \cdot)$ is non-increasing with $\lim_{t \to +\infty} \sigma(s,t) = 0$.

∗ By $\| \|_X$, we denote the norm of the normed linear space $X$. By $| |$ we denote the Euclidean norm of $\Re^n$. Let $U \subseteq X$ with $0 \in U$. By $B_U[0,r] := \{ u \in U ; \|u\|_X \leq r \}$ we denote the intersection of $U \subseteq X$ with the closed sphere of radius $r \geq 0$, centered at $0 \in U$.

∗ Let a set $U$ be a subset of a normed linear space $\mathcal{U}$, with $0 \in U$. By $\mathcal{M}(U)$ we denote the set of all locally bounded functions $u : \Re^+ \to U$. By $u_0$ we denote the identically zero input, i.e., the input that satisfies $u_0(t) = 0 \in U$ for all $t \geq 0$. If $U \subseteq \Re^n$ then $L^\infty_{loc}(\Re^+;U)$ denotes the space of measurable, locally bounded functions $u : \Re^+ \to U$.

The following convention will be adopted throughout the paper: the Cartesian product of two normed linear spaces $C := X \times Y$ will be considered to be endowed with the norm $\|(x,y)\|_C := \sqrt{\|x\|_X^2 + \|y\|_Y^2}$, unless stated otherwise. Furthermore, our results can be extended to the case of measurable and locally essentially bounded inputs (where the "sup" operator is to be understood as "essential supremum").

## 2. Input-to-Output Stability in a System-Theoretic Framework

In this section we first give the notion of a control system with outputs. We emphasize that we consider control systems which do not necessarily satisfy the classical "semigroup property" (see [16,20,38]).

**Definition 2.1:** A control system $\Sigma := (X, Y, M_U, M_D, \phi, \pi, H)$ with outputs consists of

(i) a set $U$ (control set) which is a subset of a normed linear space $\mathcal{U}$ with $0 \in U$ and a set $M_U \subseteq \mathcal{M}(U)$ (allowable control inputs) which contains at least the identically zero input $u_0$,

(ii) a set $D$ (disturbance set) and a set $M_D \subseteq \mathcal{M}(D)$, which is called the "set of allowable disturbances",

(iii) a pair of normed linear spaces $X, Y$ called the "state space" and the "output space", respectively,

(iv) a continuous map $H : \Re^+ \times X \times U \to Y$ that maps bounded sets of $\Re^+ \times X \times U$ into bounded sets of $Y$, called the "output map",

(v) a set-valued map $\Re^+ \times X \times M_U \times M_D \ni (t_0, x_0, u, d) \to \pi(t_0, x_0, u, d) \subseteq [t_0, +\infty)$, with $t_0 \in \pi(t_0, x_0, u, d)$ for all $(t_0, x_0, u, d) \in \Re^+ \times X \times M_U \times M_D$, called the set of "sampling times"

(vi) and the map $\phi : A_\phi \to X$ where $A_\phi \subseteq \Re^+ \times \Re^+ \times X \times M_U \times M_D$, called the "transition map", which has the following properties:

1) **Existence:** For each $(t_0, x_0, u, d) \in \Re^+ \times X \times M_U \times M_D$, there exists $t > t_0$ such that $[t_0, t] \times (t_0, x_0, u, d) \subseteq A_\phi$.

2) **Identity Property:** For each $(t_0, x_0, u, d) \in \Re^+ \times X \times M_U \times M_D$, it holds that $\phi(t_0, t_0, x_0, u, d) = x_0$.

3) **Causality:** For each $(t, t_0, x_0, u, d) \in A_\phi$ with $t > t_0$ and for each $(\tilde{u}, \tilde{d}) \in M_U \times M_D$ with $(\tilde{u}(\tau), \tilde{d}(\tau)) = (u(\tau), d(\tau))$ for all $\tau \in [t_0, t]$, it holds that $(t, t_0, x_0, \tilde{u}, \tilde{d}) \in A_\phi$ with $\phi(t, t_0, x_0, u, d) = \phi(t, t_0, x_0, \tilde{u}, \tilde{d})$.

4) **Weak Semigroup Property:** There exists a constant $r > 0$, such that for each $t \geq t_0$ with $(t, t_0, x_0, u, d) \in A_\phi$:

(a) $(\tau, t_0, x_0, u, d) \in A_\phi$ for all $\tau \in [t_0, t]$,
(b) $\phi(t, \tau, \phi(\tau, t_0, x_0, u, d), u, d) = \phi(t, t_0, x_0, u, d)$ for all $\tau \in [t_0, t] \cap \pi(t_0, x_0, u, d)$,
(c) if $(t + r, t_0, x_0, u, d) \in A_\phi$, then it holds that $\pi(t_0, x_0, u, d) \cap [t, t+r] \neq \emptyset$.
(d) for all $\tau \in \pi(t_0, x_0, u, d)$ with $(\tau, t_0, x_0, u, d) \in A_\phi$ we have $\pi(\tau, \phi(\tau, t_0, x_0, u, d), u, d) = \pi(t_0, x_0, u, d) \cap [\tau, +\infty)$.

In order to develop stability notions for a control system with outputs we need to clarify the notions of an equilibrium point as well as certain other important notions and classes of systems that characterize the dynamic behavior of the system (see [20,21]).

**Definition 2.2:** Let $T > 0$. A control system $\Sigma := (X, Y, M_U, M_D, \phi, \pi, H)$ with outputs is called **T-periodic**, if:

a) $H(t + T, x, u) = H(t, x, u)$ for all $(t, x, u) \in \Re^+ \times X \times U$,

b) for every $(u, d) \in M_U \times M_D$ and integer $k$ there exist inputs $P_{kT} u \in M_U$, $P_{kT} d \in M_D$ with $(P_{kT} u)(t) = u(t + kT)$ and $(P_{kT} d)(t) = d(t + kT)$ for all $t + kT \geq 0$,

c) for each $(t, t_0, x_0, u, d) \in A_\phi$ with $t \geq t_0$ and for each integer $k$ with $t_0 - kT \geq 0$ it follows that $(t - kT, t_0 - kT, x_0, P_{kT} u, P_{kT} d) \in A_\phi$ and $\pi(t_0 - kT, x_0, P_{kT} u, P_{kT} d) = \bigcup_{\tau \in \pi(t_0, x_0, u, d)} \{\tau - kT\}$ with $\phi(t, t_0, x_0, u, d) = \phi(t - kT, t_0 - kT, x_0, P_{kT} u, P_{kT} d)$.

**Definition 2.3:** A control system $\Sigma := (X, Y, M_U, M_D, \phi, \pi, H)$ with outputs is called **time-invariant or autonomous**, if:

*a)* the output map is independent of $t$, i.e., $H(t, x, u) \equiv H(x, u)$,

*b)* for every $(\theta, u, d) \in \Re \times M_U \times M_D$ there exist inputs $P_\theta u \in M_U$, $P_\theta d \in M_D$ with $(P_\theta u)(t) = u(t+\theta)$ and $(P_\theta d)(t) = d(t+\theta)$ for all $t + \theta \geq 0$,

*c)* for each $(t, t_0, x_0, u, d) \in A_\phi$ with $t \geq t_0$ and for each $\theta \in (-\infty, t_0]$ it follows that $(t-\theta, t_0 - \theta, x_0, P_\theta u, P_\theta d) \in A_\phi$ and $\pi(t_0 - \theta, x_0, P_\theta u, P_\theta d) = \bigcup_{\tau \in \pi(t_0, x_0, u, d)} \{\tau - \theta\}$ with $\phi(t, t_0, x_0, u, d) = \phi(t - \theta, t_0 - \theta, x_0, P_\theta u, P_\theta d)$.

**Definition 2.4:** Consider a control system $\Sigma := (X, Y, M_U, M_D, \phi, \pi, H)$ with outputs. We say that system $\Sigma$

(i) has the **"Boundedness-Implies-Continuation" (BIC)** property if for each $(t_0, x_0, u, d) \in \Re^+ \times X \times M_U \times M_D$, there exists a maximal existence time, i.e., there exists $t_{\max} := t_{\max}(t_0, x_0, u, d) \in (t_0, +\infty]$, such that $A_\phi = \bigcup_{(t_0, x_0, u, d) \in \Re^+ \times X \times M_U \times M_D} [t_0, t_{\max}) \times \{(t_0, x_0, u, d)\}$. In addition, if $t_{\max} < +\infty$ then for every $M > 0$ there exists $t \in [t_0, t_{\max})$ with $\|\phi(t, t_0, x_0, u, d)\|_X > M$.

(ii) is **robustly forward complete (RFC) from the input** $u \in M_U$ if it has the BIC property and for every $r \geq 0$, $T \geq 0$, it holds that

$$\sup\{\|\phi(t_0 + s, t_0, x_0, u, d)\|_X \,;\, u \in M(B_U[0, r]) \cap M_U,\, s \in [0, T],\, \|x_0\|_X \leq r,\, t_0 \in [0, T],\, d \in M_D\} < +\infty$$

**Definition 2.5:** Consider a control system $\Sigma := (X, Y, M_U, M_D, \phi, \pi, H)$ and suppose that $H(t, 0, 0) = 0$ for all $t \geq 0$. We say that $0 \in X$ is a **robust equilibrium point from the input** $u \in M_U$ for $\Sigma$ if

(i) for every $(t, t_0, d) \in \Re^+ \times \Re^+ \times M_D$ with $t \geq t_0$ it holds that $\phi(t, t_0, 0, u_0, d) = 0$.

(ii) for every $\varepsilon > 0$, $T, h \in \Re^+$ there exists $\delta := \delta(\varepsilon, T, h) > 0$ such that for all $(t_0, x, u) \in [0, T] \times X \times M_U$, $\tau \in [t_0, t_0 + h]$ with $\|x\|_X + \sup_{t \geq 0} \|u(t)\|_U < \delta$ it holds that $(\tau, t_0, x, u, d) \in A_\phi$ for all $d \in M_D$ and

$$\sup\{\|\phi(\tau, t_0, x, u, d)\|_X \,;\, d \in M_D,\, \tau \in [t_0, t_0 + h],\, t_0 \in [0, T]\} < \varepsilon$$

**Remark 2.6:** Consider a control system $\Sigma := (X, Y, M_U, M_D, \phi, \pi, H)$ with the BIC property. It follows that $\Sigma$ satisfies the (classical) semigroup property (see [20,38]) if the weak semigroup property holds with $\pi(t_0, x_0, u, d) = [t_0, t_{\max})$, where $t_{\max} \in (t_0, +\infty]$ is the maximal existence time of the transition map for $\Sigma$ that corresponds to $(t_0, x_0, u, d) \in \Re^+ \times X \times M_U \times M_D$, i.e.,

*"for each $t \in [t_0, t_{\max})$ it holds that $\phi(t, \tau, \phi(\tau, t_0, x_0, u, d), u, d) = \phi(t, t_0, x_0, u, d)$ for all $\tau \in [t_0, t]$"*
*(Classical Semigroup Property)*

The following example shows the difference between the classical semigroup property and the weak semigroup property for simple systems.

**Example 2.7:** Consider the following system:

$$\begin{aligned} \dot{x}(t) &= -x(\tau_i),\, t \in [\tau_i, \tau_{i+1}) \\ \tau_{i+1} &= \tau_i + 1 \\ x(t) &\in \Re \end{aligned} \quad (2.1)$$

with initial condition $x(t_0) = x_0 \in \Re$ and $\tau_0 = t_0 \geq 0$. Such systems will be characterized as hybrid systems with sampling partition generated by the system (see Example 2.11) and they satisfy the BIC property. In this case we can determine analytically the transition map for all $t \geq t_0$ ($u, d$ in this example are irrelevant):

$$\phi(t, t_0, x_0) = \begin{cases} (1-t+t_0)x_0 , & \text{for } t \in [t_0, t_0+1) \\ 0 , & \text{for } t \geq t_0+1 \end{cases}$$

It is clear that the state space is $\Re$ and that the classical semigroup property does not hold for this system. On the other hand the weak semigroup property holds for this system with $\pi(t_0, x_0) = \{t_0, t_0+1, t_0+2,...\}$. Notice that the set of sampling times (the sampling partition) $\pi(t_0, x_0) = \{t_0, t_0+1, t_0+2,...\}$ is generated by the system itself and depends on the initial condition. Furthermore, according to Definition 2.3, system (2.1) is autonomous.

Next, consider the following system:

$$\begin{aligned} \dot{x}(t) &= -x(\tau_i), t \in [\tau_i, \tau_{i+1}) \\ x(t) &\in \Re, \pi = \{\tau_i\}_{i=0}^{\infty} = \{0,1,2,...\} \end{aligned} \quad (2.2)$$

Such systems will be characterized as hybrid systems with impulses at fixed times (see Example 2.12) and they satisfy the BIC property. Notice that if the initial time $t_0$ is not a member of the partition $\pi = \{\tau_i\}_{i=0}^{\infty} = \{0,1,2,...\}$, then it is not possible to determine the solution of (2.2) based only on the initial condition $x(t_0) = x_0 \in \Re$ and the transition map is not well-defined. In order to be able to determine the solution of (2.2), we need to know $(x(t_0), x([t_0])) = x_0 = (x_{1,0}, x_{2,0}) \in \Re^2$ (where $[t_0]$ denotes the integer part of $t_0$). Indeed, we have:

$$x(t) = \begin{cases} x_{1,0} - (t-t_0)x_{2,0} , & t_0 \leq t < [t_0]+1 \\ (2-t+[t_0])(x_{1,0} - ([t_0]+1-t_0)x_{2,0}), & [t_0]+1 \leq t < [t_0]+2, \text{ if } t_0 \notin \pi \\ 0 , & t \geq [t_0]+2 \end{cases}$$

$$x(t) = \begin{cases} (1-t+t_0)x_{1,0} , & \text{for } t \in [t_0, t_0+1) \\ 0 , & \text{for } t \geq t_0+1 \end{cases}, \text{ if } t_0 \in \pi$$

In this case the state space is $\Re^2$ and the state of (2.2) at time $t \geq t_0$ is $\phi(t, t_0, x_0) = (x(t), x([t])) \in \Re^2$. Furthermore, notice that the classical semigroup property holds and that the partition $\pi = \{\tau_i\}_{i=0}^{\infty} = \{0,1,2,...\}$ is fixed and does not depend on the initial condition. Finally, according to Definitions 2.2 and 2.3, system (2.2) is $T$-periodic with $T = 1$ but it is not autonomous. ◁

The following examples provide classes of control systems with the BIC property and a robust equilibrium point. The examples help the reader to understand that the notions defined by Definitions 2.4-2.5 are typical for a wide class of systems under minimal assumptions.

**Example 2.8 (Finite-Dimensional Control Systems Described by Ordinary Differential Equations-ODEs):** Consider the class of systems described by ODEs of the form

$$\begin{aligned} \dot{x}(t) &= f(t, x(t), u(t), d(t)) \\ Y(t) &= H(t, x(t), u(t)) \\ x(t) &\in \Re^n, u(t) \in U, d(t) \in D, t \geq t_0 \end{aligned} \quad (2.3)$$

where $U \subseteq \Re^m$, $D \subseteq \Re^l$, with $0 \in U$ and $f : \Re^+ \times \Re^n \times U \times D \to \Re^n$, $H : \Re^+ \times \Re^n \times U \to \Re^k$ are two locally bounded mappings with $H(t,0,0) = 0$, $f(t,0,0,d) = 0$ for all $(t,d) \in \Re^+ \times D$ that satisfy the following hypotheses:

**(A1)** The mapping $(x, u, d) \to f(t, x, u, d)$ is continuous for each fixed $t \geq 0$, measurable with respect to $t \geq 0$ for each fixed $(x, u, d) \in \Re^n \times U \times D$ and such that for every bounded sets $I \subseteq \Re^+$, $S \subset \Re^n \times U$, there exists a constant $L \geq 0$ such that:

$$(x-y)'(f(t,x,u,d)-f(t,y,u,d)) \le L|x-y|^2$$
$$\forall t \in I, \forall(x,u,y,u) \in S \times S, \forall d \in D$$

**(A2)** The mapping $H : \Re^+ \times \Re^n \times U \to \Re^k$ is continuous.

**(A3)** There exist functions $\gamma \in K^+$, $a \in K_\infty$ such that $|f(t,x,u,d)| \le \gamma(t)a(|x|+|u|)$ for all $(t,x,u,d) \in \Re^+ \times \Re^n \times U \times D$.

The theory of ordinary differential equations guarantees that under hypotheses (A1-3), for each $(t_0, x_0) \in \Re^+ \times \Re^n$ and for each pair of measurable and locally bounded inputs $(u,d) \in M(U) \times M(D)$ there exists a unique absolutely continuous mapping $x(t)$ that satisfies a.e. the differential equation (2.3) with initial condition $x(t_0) = x_0 \in \Re^n$. Moreover, certain results from the theory of ordinary differential equations guarantee that (2.3) is a control system $\Sigma := (\Re^n, \Re^k, M_U, M_D, \phi, \pi, H)$ with outputs that satisfies the BIC property with $M_U, M_D$ the sets of all measurable and locally bounded mappings $u : \Re^+ \to U$, $d : \Re^+ \to D$, respectively. Furthermore, the classical semigroup property is satisfied for this system, i.e., we have $\pi(t_0, x_0, u, d) = [t_0, t_{\max})$, where $t_{\max} > t_0$ is the maximal existence time of the solution. Finally, hypotheses (A1-3) guarantee that $0 \in \Re^n$ is a robust equilibrium point from the input $u \in M_U$ for $\Sigma$. ◁

The following example presents a class of neutral functional equations described by continuous time difference equations. Such systems were recently studied in [22,35]. The importance of functional difference equations in applications is explained in [35].

**Example 2.9 (Control Systems described by Functional Difference Equations-FDEs):** Consider the class of systems described by FDEs of the form

$$x(t) = f(t, T_{r-\tau(t)}(t-\tau(t))x, u(t), d(t))$$
$$Y(t) = H(t, T_r(t)x, u(t)) \qquad (2.4)$$
$$x(t) \in \Re^n, Y(t) \in Y, u(t) \in U, d(t) \in D, t \ge t_0$$

where $r > 0$ is a constant, $\tau : \Re^+ \to (0, +\infty)$ is a positive continuous function with $\sup_{t \ge 0}\tau(t) \le r$, $D \subset \Re^l$, $U \subseteq \Re^m$ with $0 \in U$ are non-empty sets, $T_{r-\tau(t)}(t-\tau(t))x := x(t-\tau(t)+\theta)$; $\theta \in [-r+\tau(t), 0]$, $T_r(t)x := x(t+\theta)$; $\theta \in [-r, 0]$ and $H$, $f : \Omega \times U \times D \to \Re^n$, where $\Omega = \cup_{t \ge 0}\{t\} \times \mathcal{F}_t$ and $\mathcal{F}_t$ denotes the set of bounded functions $x : [-r+\tau(t), 0] \to \Re^n$, are locally bounded mappings which satisfy the following hypotheses:

**(R1)** There exist functions $\gamma \in K^+$, $a \in K_\infty$ such that $\left|f(t, T_{r-\tau(t)}(-\tau(t))x, u, d)\right| \le \gamma(t)a\left(\sup_{\theta \in [-r, -\tau(t)]}|x(\theta)| + |u|\right)$ for all $(t, x, u, d) \in \Re^+ \times X \times U \times D$, where $X$ is the normed linear space of the bounded functions $x : [-r, 0] \to \Re^n$ with $\|x\|_X := \sup_{\theta \in [-r, 0]}|x(\theta)|$.

**(R2)** The output map $H : \Re^+ \times X \times U \to Y$, where $Y$ is a normed linear space, is a continuous mapping that maps bounded sets of $\Re^+ \times X \times U$ into bounded sets of $Y$ with $H(t, 0, 0) = 0$ for all $t \ge 0$.

It should be clear that under the hypotheses stated above, for each $(t_0, x_0) \in \Re^+ \times X$ and for each pair of locally bounded functions $u : \Re^+ \to U$, $d : \Re^+ \to D$ there exists a unique locally bounded mapping $x(t)$ that satisfies the difference equations (2.4) with initial condition $x(t_0 + \theta) = x_0(\theta); \theta \in [-r, 0]$. Consequently, (2.4) describes a control system $\Sigma := (X, Y, M_U, M_D, \phi, \pi, H)$ with outputs and evolution map $\phi$ defined by $\phi(t, t_0, x_0, u, d) = x(t+\theta)$; $\theta \in [-r, 0]$, where $\mathcal{U} := \mathfrak{R}^m$, $M_U$ the set of all locally bounded functions $u : \mathfrak{R}^+ \to U$ and $M_D$ the set of all functions $d : \mathfrak{R}^+ \to D$.

Systems described by functional difference evolution equations of the form (2.4) are considered in [4,22,35]. Working exactly in the same way as in [22], it can be shown that system (2.4) is RFC from the input $u \in M_U$ and that $0 \in \mathcal{X}$ is a robust equilibrium point from the input $u \in M_U$ for system (2.4).

Notice that a major advantage of allowing the output to take values in abstract normed linear spaces is that we are in a position to consider:

- outputs with no delays, e.g. $Y(t) = h(t, x(t), u(t))$ with $Y = \mathfrak{R}^k$,
- outputs with discrete or distributed delay, e.g. $Y(t) = h(t, x(t), x(t-r), u(t))$ or $Y(t) = \sup_{\theta \in [t-r, t]} h(t, \theta, x(\theta), u(t))$ with $Y = \mathfrak{R}^k$,
- functional outputs with memory, e.g. $Y(t) = h(t, \theta, x(t+\theta))$; $\theta \in [-r, 0]$ or the identity output $Y(t) = T_r(t)x = x(t+\theta)$; $\theta \in [-r, 0]$ with $Y = \mathcal{X}$.

Finally, notice that the classical semigroup property is satisfied for this system, i.e., we have $\pi(t_0, x_0, u, d) = [t_0, +\infty)$. ◁

The following example is an immediate consequence of Theorems 2.2 and 3.2 in [4], concerning continuous dependence on initial conditions and continuation of solutions of retarded functional differential equations, respectively.

**Example 2.10 (Control Systems described by Retarded Functional Differential Equations-RFDEs):** Consider the class of systems described by RFDEs of the form

$$\dot{x}(t) = f(t, T_r(t)x, u(t), d(t))$$
$$Y(t) = H(t, T_r(t)x, u(t)) \quad (2.5)$$
$$x(t) \in \mathfrak{R}^n, \ u(t) \in U, \ d(t) \in D, \ t \geq t_0$$

where $T_r(t)x := x(t+\theta)$; $\theta \in [-r, 0]$, $D \subseteq \mathfrak{R}^l$ is a non-empty set, $U \subseteq \mathfrak{R}^m$ is a non-empty set with $0 \in U$, $f : \mathfrak{R}^+ \times C^0([-r,0]; \mathfrak{R}^n) \times U \times D \to \mathfrak{R}^n$, $H : \mathfrak{R}^+ \times C^0([-r,0]; \mathfrak{R}^n) \times U \to Y$ ($Y$ is a normed linear space) are locally bounded mappings with $f(t, 0, 0, d) = 0$, $H(t, 0, 0) = 0$ for all $(t, d) \in \mathfrak{R}^+ \times D$, that satisfy the following hypotheses:

**(S1)** The mapping $(x, u, d) \to f(t, x, u, d)$ is continuous for each fixed $t \geq 0$ and such that for every bounded $I \subseteq \mathfrak{R}^+$ and for every bounded $S \subset C^0([-r,0]; \mathfrak{R}^n) \times U$, there exists a constant $L \geq 0$ such that:

$$(x(0) - y(0))'(f(t, x, u, d) - f(t, y, u, d)) \leq L \max_{\tau \in [-r, 0]} |x(\tau) - y(\tau)|^2$$
$$\forall t \in I, \ \forall (x, u, y, u) \in S \times S, \ \forall d \in D$$

**(S2)** There exist functions $\gamma \in K^+$, $a \in K_\infty$ such that $|f(t, x, u, d)| \leq \gamma(t) a(\|x\|_r + |u|)$ for all $(t, x, u, d) \in \mathfrak{R}^+ \times C^0([-r,0]; \mathfrak{R}^n) \times U \times D$, where $\|x\|_r$ denotes the sup-norm of the space $C^0([-r,0]; \mathfrak{R}^n)$, i.e., $\|x\|_r := \max_{\theta \in [-r, 0]} |x(\theta)|$.

**(S3)** There exists a countable set $A \subset \mathfrak{R}^+$, which is either finite or $A = \{t_k \ ; k = 1, ..., \infty\}$ with $t_{k+1} > t_k > 0$ for all $k = 1, 2, ...$ and $\lim t_k = +\infty$, such that mapping $(t, x, u, d) \in (\mathfrak{R}^+ \setminus A) \times C^0([-r,0]; \mathfrak{R}^n) \times U \times D \to f(t, x, u, d)$ is continuous. Moreover, for each fixed $(t_0, x, u, d) \in \Re^+ \times C^0([-r,0]; \Re^n) \times U \times D$, we have $\lim_{t \to t_0^+} f(t, x, u, d) = f(t_0, x, u, d)$.

**(S4)** The mapping $H : \Re^+ \times C^0([-r,0]; \Re^n) \times U \to Y$ is a continuous mapping that maps bounded sets of $\Re^+ \times C^0([-r,0]; \Re^n) \times U$ into bounded sets of $Y$.

The theory of retarded functional differential equations guarantees that under hypotheses (S1-4), for each $(t_0, x_0) \in \Re^+ \times C^0([-r,0]; \Re^n)$ and for each pair of measurable and locally bounded inputs $(u,d) \in \mathcal{M}(U) \times \mathcal{M}(D)$ there exists a unique absolutely continuous mapping $x(t)$ that satisfies a.e. the differential equation (2.5) with initial condition $x(t_0) = x_0 \in C^0([-r,0]; \Re^n)$. Moreover, certain results from the theory of retarded functional differential equations (Theorems 2.2 and 3.2 in [4]) guarantee that (2.5) is a control system $\Sigma := (C^0([-r,0]; \Re^n), Y, M_U, M_D, \phi, \pi, H)$ with outputs that satisfies the BIC property with $M_U, M_D$ the sets of all measurable and locally bounded mappings $u : \Re^+ \to U$, $d : \Re^+ \to D$, respectively. Furthermore, the classical semigroup property is satisfied for this system, i.e., we have $\pi(t_0, x_0, u, d) = [t_0, t_{\max})$, where $t_{\max} > t_0$ is the maximal existence time of the solution. Finally, hypotheses (S1-4) guarantee that $0 \in C^0([-r,0]; \Re^n)$ is a robust equilibrium point from the input $u \in M_U$ for $\Sigma$. Again notice that a major advantage of allowing the output to take values in abstract normed linear spaces is that we are in a position to consider various output cases (see previous example). ◁

The following example presents an important class of systems that does not satisfy the classical "semigroup property".

**Example 2.11 (Hybrid Systems with sampling partition generated by the system):** Consider the class of systems described by impulsive differential equations of the form

$$\dot{x}(t) = f(t, \tau_i, x(t), x(\tau_i), u(t), u(\tau_i), d(t), d(\tau_i)) \quad , \quad t \in [\tau_i, \tau_{i+1})$$
$$\tau_0 = t_0, \tau_{i+1} = \tau_i + h(\tau_i, x(\tau_i), u(\tau_i), d(\tau_i)), i = 0,1,...$$
$$x(\tau_{i+1}) = R\left(\tau_i, \lim_{t \to \tau_{i+1}^-} x(t), x(\tau_i), u(\tau_{i+1}), u(\tau_i), d(\tau_{i+1}), d(\tau_i)\right) \quad (2.6)$$
$$Y(t) = H(t, x(t), u(t))$$

where $D \subseteq \Re^l$, $U \subseteq \Re^m$ is a closed set with $0 \in U$, $h : \Re^+ \times \Re^n \times U \times D \to (0, r]$ is a positive function which is bounded by certain constant $r > 0$, $f : \Re^+ \times \Re^+ \times \Re^n \times \Re^n \times U \times U \times D \times D \to \Re^n$, $H : \Re^+ \times \Re^n \times U \to \Re^p$ and $R : \Re^+ \times \Re^n \times \Re^n \times U \times U \times D \times D \to \Re^n$ is a triplet of vector fields that satisfy the following hypotheses:

**(P1)** $f(t, \tau, x, x_0, u, u_0, d, d_0)$ is measurable with respect to $t \geq 0$, continuous with respect to $(x, d, u) \in \Re^n \times D \times U$ and such that for every bounded $S \subset \Re^+ \times \Re^+ \times \Re^n \times \Re^n \times U \times U$ there exists constant $L \geq 0$ such that

$$(x-y)'\bigl(f(t, \tau, x, x_0, u, u_0, d, d_0) - f(t, \tau, y, x_0, u, u_0, d, d_0)\bigr) \leq L|x-y|^2$$
$$\forall (t, \tau, x, x_0, u, u_0, d, d_0) \in S \times D \times D, \forall (t, \tau, y, x_0, u, u_0, d, d_0) \in S \times D \times D$$

**(P2)** There exist functions $\gamma \in K^+$, $a \in K_\infty$ such that

$$|f(t, \tau, x, x_0, u, u_0, d, d_0)| \leq \gamma(t) a(|x| + |x_0| + |u| + |u_0|), \forall (\tau, u, u_0, d, d_0, x, x_0) \in \Re^+ \times U \times U \times D \times D \times \Re^n \times \Re^n, \forall t \geq \tau$$

$$|R(t, x, x_0, u, u_0, d, d_0)| \leq \gamma(t) a(|x| + |x_0| + |u| + |u_0|), \forall (t, u, u_0, d, d_0, x, x_0) \in \Re^+ \times U \times U \times D \times D \times \Re^n \times \Re^n$$

**(P3)** $H : \Re^+ \times \Re^n \times U \to \Re^p$ is a continuous map with $H(t, 0, 0) = 0$ for all $t \geq 0$,

**(P4)** There exist a positive, continuous and bounded function $h_l : \Re^+ \times \Re^n \times U \to (0, r]$ and a partition $\pi = \{T_i\}_{i=0}^{\infty}$ of $\Re^+$, i.e., an increasing sequence of times with $T_0 = 0$ and $T_i \to +\infty$ such that:

$$h(t, x, u, d) \geq \min\{p_\pi(t) - t, h_l(t, x, u)\}, \quad \forall (t, x, u, d) \in \Re^+ \times \Re^n \times U \times D$$

where $p_\pi(t) := \min\{T \in \pi ; t < T\}$.

Hybrid systems of the form (2.6) under hypotheses (P1-4) are considered in [21,22], where it is shown that for each $(t_0, x_0) \in \Re^+ \times \Re^n$ and for each pair of measurable and locally bounded inputs $u : \Re^+ \to U$ and $d : \Re^+ \to D$ there exists a unique piecewise absolutely continuous function $t \to x(t) \in \Re^n$ with initial condition $x(t_0) = x_0$, which is produced by the following algorithm:

Step $i$:
1) Given $\tau_i$ and $x(\tau_i)$, calculate $\tau_{i+1}$ using the equation $\tau_{i+1} = \tau_i + h(\tau_i, x(\tau_i), u(\tau_i), d(\tau_i))$,
2) Compute the state trajectory $x(t)$, $t \in [\tau_i, \tau_{i+1})$ as the solution of the differential equation $\dot{x}(t) = f(t, \tau_i, x(t), x(\tau_i), u(t), u(\tau_i), d(t), d(\tau_i))$,
3) Calculate $x(\tau_{i+1})$ using the equation $x(\tau_{i+1}) = R\left(\tau_i, \lim_{t \to \tau_{i+1}^-} x(t), x(\tau_i), u(\tau_{i+1}), u(\tau_i), d(\tau_{i+1}), d(\tau_i)\right)$,
4) Compute the output trajectory $Y(t)$, $t \in [\tau_i, \tau_{i+1}]$ using the equation $Y(t) = H(t, x(t), u(t))$

For $i = 0$ we take $\tau_0 = t_0$ and $x(\tau_0) = x_0$ (initial condition).

In [21] it is shown that system (2.6) under hypotheses (P1-4) is a control system $\Sigma := (\mathcal{X}, Y, M_U, M_D, \phi, \pi, H)$ with outputs with the BIC property for which $0 \in \Re^n$ is a robust equilibrium point from the input $u \in M_U$. Particularly, we have $\mathcal{X} = \Re^n$, $Y = \Re^p$, $\mathcal{U} = \Re^m$ and $M_U$, $M_D$ the sets of measurable and locally bounded inputs $u : \Re^+ \to U$ and $d : \Re^+ \to D$, respectively. The set $\pi(t_0, x_0, u, d) \subseteq [t_0, +\infty)$ involved in the weak semigroup property consists of the sequence $\pi = \{\tau_i\}_{i=0}^{\infty}$ generated by the recursive relation $\tau_{i+1} = \tau_i + h(\tau_i, x(\tau_i), u(\tau_i), d(\tau_i))$, $i = 0, 1, \ldots$ with $\tau_0 = t_0$. Notice that the control system (2.6) fails to satisfy the classical semigroup property.

Notice that if $h(\tau + T, x, u, d) = h(\tau, x, u, d)$, $f(t + T, \tau + T, x, x_0, u, u_0, d, d_0) = f(t, \tau, x, x_0, u, u_0, d, d_0)$, $R(\tau + T, x, x_0, u, u_0, d, d_0) = R(\tau, x, x_0, u, u_0, d, d_0)$ and $H(t + T, x, u) = H(t, x, u)$ for certain $T > 0$ and for $(t, \tau, u, u_0, d, d_0, x, x_0) \in \Re^+ \times \Re^+ \times U \times U \times D \times D \times \Re^n \times \Re^n$ with $t \geq \tau$, then system (2.6) is $T$-periodic. Moreover, if $h(\tau, x, u, d) = h(x, u, d)$, $f(t, \tau, x, x_0, u, u_0, d, d_0) = f(t - \tau, x, x_0, u, u_0, d, d_0)$, $R(\tau, x, x_0, u, u_0, d, d_0) = R(x, x_0, u, u_0, d, d_0)$ and $H(t, x, u) = H(x, u)$ for $(t, \tau, u, u_0, d, d_0, x, x_0) \in \Re^+ \times \Re^+ \times U \times U \times D \times D \times \Re^n \times \Re^n$ with $t \geq \tau$ then system (2.6) is autonomous.

Systems of the form (2.6) under hypotheses (P1-4) arise frequently in certain applications in mathematical control theory and numerical analysis. Specifically, they arise when

i) a (not necessarily continuous) sampled-data feedback law (with possibly variable sampling rate) is applied to a finite-dimensional control system. For example, state-dependent sampling rates were related in [2] with the classical work on discontinuous stabilizability in [1] while feedback stabilization problems with zero order hold and constant positive sampling rate were considered in [28-34] and time-varying sampling rates were considered in [6,7].
ii) a numerical discretization method (with possibly variable integration step sizes) is applied in order to obtain the numerical solution of a given system of ordinary differential equations; see [3,41] for the case of constant integration step sizes and [19,21] for the case of variable integration step sizes.

For a unified description of the above problems, see [21,22].  ◁

In contrast with the previous example, it should be noted that hybrid systems with impulses at fixed times satisfy the classical "semigroup property". The following example illustrates this case.

**Example 2.12 (Hybrid Systems with Impulses at fixed times):** Consider the class of systems described by impulsive differential equations of the form

$$\dot{x}(t) = f(t, d(t), d(\tau_i), x(t), x(\tau_i), u(t), u(\tau_i)) \quad , \quad \tau_i \le t < \tau_{i+1}$$

$$x(\tau_{i+1}) = R\left(\tau_i, \lim_{t \to \tau_{i+1}^-} x(t), x(\tau_i), u(\tau_{i+1}), u(\tau_i), d(\tau_{i+1}), d(\tau_i)\right) \quad (2.7)$$

$$Y(t) = H(t, x(t), u(t))$$

$$x(t) \in \Re^n \;,\; Y(t) \in \Re^k \;,\; u(t) \in V \subseteq \Re^m \;,\; t \ge 0 \;,\; d(t) \in D$$

where $D \subseteq \Re^l$, $V \subseteq \Re^m$ is a closed set with $0 \in V$, $\pi = \{\tau_i\}_{i=0}^\infty$ is a partition of $\Re^+$ with diameter $r > 0$, i.e., an increasing sequence of times with $\tau_0 = 0$, $\sup\{\tau_{i+1} - \tau_i \;;\; i = 0,1,2,...\} = r$ and $\tau_i \to +\infty$, $d(t)$ represents the disturbance vector or the vector of time-varying uncertainties taking values in the set $D \subset \Re^l$, $Y(t)$ represents the output of the system and $u(t) \in V$ represents the input vector. A wide class of systems described by impulsive differential equations with impulses at fixed times, as well as hybrid systems of the form:

$$\dot{x}(t) = f(t, x(t), u(t), w(i)) \quad , \quad \tau_i \le t < \tau_{i+1}$$

$$w(i) = g(i, x(\tau_i), u(\tau_i)) \quad (2.8)$$

where $\pi = \{\tau_i\}_{i=0}^\infty$ is a partition of $\Re^+$ of diameter $r > 0$, can be represented by the time-varying case (2.7). Fundamental properties of the solutions of systems of the form (2.8) are studied in [24,25].

Consider system (2.7) under the following assumptions:

**(Q1)** $\pi = \{\tau_i\}_{i=0}^\infty$ is a partition of $\Re^+$ with finite diameter $r > 0$, i.e., an increasing sequence of times with $\tau_0 = 0$, $\sup\{\tau_{i+1} - \tau_i \;;\; i = 0,1,2,...\} = r$ and $\tau_i \to +\infty$.

**(Q2)** $H : \Re^+ \times \Re^n \times V \to \Re^k$ is continuous with $H(t,0,0) = 0$, for all $t \ge 0$.

**(Q3)** $f(t, d, d_0, x, x_0, u, u_0)$ is measurable with respect to $t \ge 0$, continuous with respect to $(x, d, u) \in \Re^n \times D \times V$ and such that for every compact $S \subset \Re^n \times \Re^n \times V \times V$ and for every compact $I \subset \Re^+$ there exists constant $L \ge 0$ such that

$$(x - y)'(f(t, d, d_0, x, x_0, u, u_0) - f(t, d, d_0, y, x_0, u, u_0)) \le L|x - y|^2$$

$$\forall t \in I \;,\; \forall (d, d_0) \in D \times D \;,\; \forall (x, x_0, u, u_0) \in S \;,\; \forall (y, x_0, u, u_0) \in S$$

**(Q4)** There exist functions $\gamma \in K^+$, $a \in K_\infty$ such that $|f(t, d, d_0, x, x_0, u, u_0)| \le \gamma(t) a(|x_0| + |x| + |u| + |u_0|)$, $|R(t, x, x_0, u, u_0, d, d_0)| \le \gamma(t) a(|x| + |x_0| + |u| + |u_0|)$, for all $(t, d, d_0, x, x_0, u, u_0) \in \Re^+ \times D \times D \times \Re^n \times \Re^n \times V \times V$.

Systems of the form (2.7) with $R(t, x, x_0, u, u_0, d, d_0) \equiv x$ (impulse free case) were considered in [23]. Special classes of impulsive systems of the form (2.7) were studied in [5]. Using the method of steps on consecutive intervals, it is clear that system (2.7) under hypotheses (Q1-4) defines a control system $\Sigma := (\mathcal{X}, Y, M_U, M_D, \phi, \pi, H)$ with outputs and the BIC property, with state space $\mathcal{X} = \Re^n \times \Re^n$, output space $Y = \Re^k$, set of structured uncertainties $M_D$ being the set of mappings $t \in \Re^+ \to d(t) = \{\tilde{d}(t + \theta) ; \theta \in [-r, 0]\}$ where $\tilde{d} : \Re \to D$ is any measurable and locally bounded function, input space $\mathcal{U}$ the normed linear space of measurable and bounded functions on $[-r, 0]$ taking values in $\Re^m$ endowed with the sup norm, $U \subseteq \mathcal{U}$ the set of measurable and bounded functions on $[-r, 0]$ taking values in $V \subseteq \Re^m$ and set of external inputs $M_U$ being the set of mappings $t \in \Re^+ \to u(t) = \{\tilde{u}(t+\theta) ; \theta \in [-r,0]\} \in U$, where $\tilde{u} : \Re \to \Re^m$ is a measurable and locally bounded function. The reader may be surprised by the complicated definition of $M_D$ and $M_U$, but it should be emphasized that this definition guarantees that the causality property of the control system (2.7) holds. Notice that the classical semigroup property is satisfied for this system, i.e., we have $\pi(t_0, x_0, u, d) = [t_0, t_{max})$, where $t_{max} > t_0$ is the maximal existence time of the solution. However notice that if the vector fields $f$ and $R$ are independent of $d(\tau_i), x(\tau_i), u(\tau_i)$ (this is the case studied in [5]) then system (2.7) under hypotheses (Q1-4) defines a control system $\Sigma := (X, Y, M_U, M_D, \phi, \pi, H)$ with outputs and the BIC property, with state space $X = \Re^n$, output space $Y = \Re^k$, set of structured uncertainties $M_D$ being the set of measurable and locally bounded functions $d : \Re \to D$, input space $U = \Re^m$ and $M_U$ being the set of measurable and locally bounded functions $u : \Re \to U$.

Let $q_\pi(t) = \max\{\tau_i ; \tau_i \in \pi, \tau_i \leq t\}$. For all $(t_0, x_0, x_1, d, u) \in \Re^+ \times \Re^n \times \Re^n \times M_D \times M_U$, we denote by $x(t) = \phi(t, t_0, x_0, x_1; d, u) \in \Re^n$ the solution of (2.7) at time $t \geq t_0$ with initial condition $x(t_0) = x_0$ and the additional condition $x(q_\pi(t_0)) = x_1$, which holds only for the case $t_0 \notin \pi$, corresponding to inputs $(d, u) \in M_D \times M_U$ (this solution is unique by virtue of property (Q3)). Notice that the actual state of system (2.7) at time $t \geq t_0$ is given by $\tilde{\phi}(t, t_0, x_0, x_1; d, u) = (\phi(t, t_0, x_0, x_1; d, u), \phi(q_\pi(t), t_0, x_0, x_1; d, u)) \in \Re^n \times \Re^n$.

Hypotheses (Q3-4) can be used in order to show that $0 \in \Re^n \times \Re^n$ is a robust equilibrium point from the input $u \in M_U$, exactly in the same way with the proof of the analogous result in [23]. Notice that if $f(t+T, x, x_0, u, u_0, d, d_0) = f(t, x, x_0, u, u_0, d, d_0)$, $R(t+T, x, x_0, u, u_0, d, d_0) = R(t, x, x_0, u, u_0, d, d_0)$, $H(t+T, x, u) = H(t, x, u)$ and $\pi = \{iT\}_{i=0}^\infty$ for certain $T > 0$ and for all $(t, u, u_0, d, d_0, x, x_0) \in \Re^+ \times V \times V \times D \times D \times \Re^n \times \Re^n$, then system (2.7) is $T$-periodic. Moreover, it should be noted that there system (2.7) fails to be autonomous for every possible selection of the sets $D$, $V$, vector fields $f, R, H$ and partition $\pi$. ◁

For control systems with the BIC property the following lemma provides a useful characterization of the RFC property. Its proof is provided in the Appendix.

**Lemma 2.13:** *System $\Sigma := (X, Y, M_U, M_D, \phi, \pi, H)$ is RFC from the input $u \in M_U$ if and only if for every $\beta \in K^+$ there exist functions $\mu, c \in K^+$, $a, p \in K_\infty$ (depending only on $\beta \in K^+$) such that the following estimate holds for all $(t_0, x_0, d, u) \in \Re^+ \times X \times M_D \times M_U$:*

$$\beta(t)\|\phi(t, t_0, x_0, u, d)\|_X \leq \max\left\{\mu(t-t_0), c(t_0), a(\|x_0\|_X), \sup_{t_0 \leq \tau \leq t} p(\|u(\tau)\|_U)\right\}, \forall t \geq t_0 \qquad (2.9)$$

Next we present the Input-to-Output Stability property for the class of systems described by Definition 2.1.

**Definition 2.14:** *Consider a control system $\Sigma := (X, Y, M_U, M_D, \phi, \pi, H)$ with outputs and the BIC property and for which $0 \in X$ is a robust equilibrium point from the input $u \in M_U$. Suppose that $\Sigma$ is RFC from the input $u \in M_U$.*

\* *If there exist functions $\sigma \in KL$, $\beta, \delta \in K^+$, $\gamma \in N$ such that the following estimate holds for all $u \in M_U$, $(t_0, x_0, d) \in \Re^+ \times X \times M_D$ and $t \geq t_0$:*

$$\|H(t, \phi(t, t_0, x_0, u, d), u(t))\|_Y \leq \sigma(\beta(t_0)\|x_0\|_X, t-t_0) + \sup_{t_0 \leq \tau \leq t} \gamma(\delta(\tau)\|u(\tau)\|_U) \qquad (2.10)$$

*then we say that $\Sigma$ satisfies the Weighted Input-to-Output Stability (WIOS) property from the input $u \in M_U$ with gain $\gamma \in N$ and weight $\delta \in K^+$. Moreover, if $\beta(t) \equiv 1$ then we say that $\Sigma$ satisfies the Uniform Weighted Input-to-Output Stability (UWIOS) property from the input $u \in M_U$ with gain $\gamma \in N$ and weight $\delta \in K^+$.*

∗ *If there exist functions* $\sigma \in KL$, $\beta \in K^+$, $\gamma \in N$ *such that the following estimate holds for all* $u \in M_U$, $(t_0, x_0, d) \in \Re^+ \times X \times M_D$ *and* $t \geq t_0$:

$$\|H(t, \phi(t, t_0, x_0, u, d), u(t))\|_Y \leq \sigma(\beta(t_0)\|x_0\|_X, t - t_0) + \sup_{t_0 \leq \tau \leq t} \gamma(\|u(\tau)\|_U) \quad (2.11)$$

*then we say that* $\Sigma$ *satisfies the Input-to-Output Stability (IOS) property from the input* $u \in M_U$ *with gain* $\gamma \in N$. *Moreover, if* $\beta(t) \equiv 1$ *then we say that* $\Sigma$ *satisfies the Uniform Input-to-Output Stability (UIOS) property from the input* $u \in M_U$ *with gain* $\gamma \in N$.

*Finally, for the special case of the identity output mapping, i.e.,* $H(t, x, u) := x$, *the (Uniform) (Weighted) Input-to-Output Stability property from the input* $u \in M_U$ *is called (Uniform) (Weighted) Input-to-State Stability property from the input* $u \in M_U$.

**Remark 2.15:** Using the inequalities $\max\{a, b\} \leq a + b \leq \max\{a + \rho(a), b + \rho^{-1}(b)\}$ (which hold for all $\rho \in K_\infty$ and $a, b \geq 0$), it should be clear that the WIOS property for $\Sigma := (X, Y, M_U, M_D, \phi, \pi, H)$ can be defined by using an estimate of the form:

$$\|H(t, \phi(t, t_0, x_0, u, d), u(t))\|_Y \leq \max\left\{\sigma(\beta(t_0)\|x_0\|_X, t - t_0), \sup_{t_0 \leq \tau \leq t} \gamma(\delta(\tau)\|u(\tau)\|_U)\right\} \quad (2.10')$$

instead of (2.10). Similarly, the IOS property for $\Sigma := (X, Y, M_U, M_D, \phi, \pi, H)$ can be defined by using an estimate of the form:

$$\|H(t, \phi(t, t_0, x_0, u, d), u(t))\|_Y \leq \max\left\{\sigma(\beta(t_0)\|x_0\|_X, t - t_0), \sup_{t_0 \leq \tau \leq t} \gamma(\|u(\tau)\|_U)\right\} \quad (2.11')$$

instead of (2.11).

The following lemmas provide $\varepsilon - \delta$ characterizations of the WIOS and UWIOS properties, which are going to be used in the following section of the paper. Their proofs are provided in the Appendix.

**Lemma 2.16:** *Consider a control system* $\Sigma := (X, Y, M_U, M_D, \phi, \pi, H)$ *with outputs and the BIC property and for which* $0 \in X$ *is a robust equilibrium point from the input* $u \in M_U$. *Suppose that* $\Sigma$ *is RFC from the input* $u \in M_U$. *Furthermore, suppose that there exist functions* $V : \Re^+ \times X \times U \to \Re^+$ *with* $V(t, 0, 0) = 0$ *for all* $t \geq 0$, $\gamma \in N$ *and* $\delta \in K^+$ *such that the following properties hold:*

**P1** *For every* $s \geq 0$, $T \geq 0$, *it holds that*

$$\sup\left\{V(t, \phi(t, t_0, x_0, u, d), u(t)) - \sup_{t_0 \leq \tau \leq t} \gamma(\delta(\tau)\|u(\tau)\|_U); t \geq t_0, \|x_0\|_X \leq s, t_0 \in [0, T], d \in M_D, u \in M_U\right\} < +\infty$$

**P2** *For every* $\varepsilon > 0$ *and* $T \geq 0$ *there exists a* $\delta := \delta(\varepsilon, T) > 0$, *such that:*

$$\sup\left\{V(t, \phi(t, t_0, x_0, u, d), u(t)) - \sup_{t_0 \leq \tau \leq t} \gamma(\delta(\tau)\|u(\tau)\|_U); t \geq t_0, \|x_0\|_X \leq \delta, t_0 \in [0, T], d \in M_D, u \in M_U\right\} \leq \varepsilon$$

**P3** *For every* $\varepsilon > 0$, $T \geq 0$ *and* $R \geq 0$, *there exists a* $\tau := \tau(\varepsilon, T, R) \geq 0$, *such that:*

$$\sup\left\{V(t, \phi(t, t_0, x_0, u, d), u(t)) - \sup_{t_0 \leq \tau \leq t} \gamma(\delta(\tau)\|u(\tau)\|_U); t \geq t_0 + \tau, \|x_0\|_X \leq R, t_0 \in [0, T], d \in M_D, u \in M_U\right\} \leq \varepsilon$$

*Then there exist functions* $\sigma \in KL$ *and* $\beta \in K^+$ *such that the following estimate holds for all* $u \in M_U$, $(t_0, x_0, d) \in \Re^+ \times X \times M_D$ *and* $t \geq t_0$:

$$V(t,\phi(t,t_0,x_0,u,d),u(t)) \le \sigma\big(\beta(t_0)\|x_0\|_\chi, t-t_0\big) + \sup_{t_0 \le \tau \le t} \gamma\big(\delta(\tau)\|u(\tau)\|_\mathcal{U}\big) \quad (2.12)$$

*Moreover, if there exists* $a \in N$ *such that* $\|H(t,x,u)\|_Y \le a(V(t,x,u))$ *for all* $(t,x,u) \in \Re^+ \times \chi \times U$, *then for every* $\rho \in K_\infty$, $\Sigma$ *satisfies the WIOS property from the input* $u \in M_U$ *with gain* $\tilde\gamma \in N$ *and weight* $\delta \in K^+$, *where* $\tilde\gamma(s) := a(\gamma(s) + \rho(\gamma(s)))$.

**Lemma 2.17:** *Consider a control system* $\Sigma := (\chi, Y, M_U, M_D, \phi, \pi, H)$ *with outputs and the BIC property and for which* $0 \in \chi$ *is a robust equilibrium point from the input* $u \in M_U$. *Suppose that* $\Sigma$ *is RFC from the input* $u \in M_U$. *Furthermore, suppose that there exist functions* $V : \Re^+ \times \chi \times U \to \Re^+$ *with* $V(t,0,0) = 0$ *for all* $t \ge 0$, $\gamma \in N$ *and* $\delta \in K^+$ *such that the following properties hold:*

**P1** *For every* $s \ge 0$, *it holds that*

$$\sup\left\{ V(t,\phi(t,t_0,x_0,u,d),u(t)) - \sup_{t_0 \le \tau \le t} \gamma\big(\delta(\tau)\|u(\tau)\|_\mathcal{U}\big); t \ge t_0, \|x_0\|_\chi \le s, t_0 \ge 0, d \in M_D, u \in M_U \right\} < +\infty$$

**P2** *For every* $\varepsilon > 0$ *there exists a* $\delta := \delta(\varepsilon) > 0$, *such that:*

$$\sup\left\{ V(t,\phi(t,t_0,x_0,u,d),u(t)) - \sup_{t_0 \le \tau \le t} \gamma\big(\delta(\tau)\|u(\tau)\|_\mathcal{U}\big); t \ge t_0, \|x_0\|_\chi \le \delta, t_0 \ge 0, d \in M_D, u \in M_U \right\} \le \varepsilon$$

**P3** *For every* $\varepsilon > 0$ *and* $R \ge 0$, *there exists a* $\tau := \tau(\varepsilon, R) \ge 0$, *such that:*

$$\sup\left\{ V(t,\phi(t,t_0,x_0,u,d),u(t)) - \sup_{t_0 \le \tau \le t} \gamma\big(\delta(\tau)\|u(\tau)\|_\mathcal{U}\big); t \ge t_0 + \tau, \|x_0\|_\chi \le R, t_0 \ge 0, d \in M_D, u \in M_U \right\} \le \varepsilon$$

*Then there exists a function* $\sigma \in KL$ *such that estimate (2.12) holds for all* $u \in M_U$, $(t_0, x_0, d) \in \Re^+ \times \chi \times M_D$ *and* $t \ge t_0$ *with* $\beta(t) \equiv 1$. *Moreover, if there exists* $a \in N$ *such that* $\|H(t,x,u)\|_Y \le a(V(t,x,u))$ *for all* $(t,x,u) \in \Re^+ \times \chi \times U$, *then for every* $\rho \in K_\infty$, $\Sigma$ *satisfies the UWIOS property from the input* $u \in M_U$ *with gain* $\tilde\gamma \in N$ *and weight* $\delta \in K^+$, *where* $\tilde\gamma(s) := a(\gamma(s) + \rho(\gamma(s)))$.

**Remark:** Notice that Lemmata 2.16 and 2.17 can be very useful for the demonstration of the (U)WIOS property, because in practice we show properties (P1-3) for some Lyapunov function $V$ and not necessarily for the norm of the output map. Moreover, notice that it is not required that $V$ is continuous. If $V : \Re^+ \times \chi \times U \to \Re^+$ is a continuous functional that maps bounded sets of $\Re^+ \times \chi \times U$ into bounded sets of $\Re^+$, then Lemmata 2.16 and 2.17 guarantee that $\Sigma$ satisfies the WIOS and the UWIOS properties with $V$ as output, respectively, from the input $u \in M_U$ with gain $\gamma \in N$ and weight $\delta \in K^+$.

Finally, we end this section with some useful observations for *T-periodic* control systems. It turns out that periodicity guarantees uniformity with respect to the initial times. The following lemmas should be compared with Lemma 1.1, page 131 in [4]. Their proofs are provided in the Appendix.

**Lemma 2.18:** *Suppose that* $\Sigma := (\chi, Y, M_U, M_D, \phi, \pi, H)$ *is T-periodic. If* $\Sigma$ *satisfies the WIOS property from the input* $u \in M_U$, *then* $\Sigma$ *satisfies the UWIOS property from the input* $u \in M_U$.

**Lemma 2.19:** *Suppose that* $\Sigma := (\chi, Y, M_U, M_D, \phi, \pi, H)$ *is T-periodic. If* $\Sigma$ *satisfies the IOS property from the input* $u \in M_U$, *then* $\Sigma$ *satisfies the UIOS property from the input* $u \in M_U$.

## 3. A Small-Gain Theorem for a Wide Class of Systems

The main result of the present work is stated next.

**Theorem 3.1:** *Consider the system* $\Sigma := (\mathcal{X}, \mathcal{Y}, M_U, M_D, \phi, \pi, H)$ *with the BIC property and suppose that there exist maps* $V_1 : \Re^+ \times \mathcal{X} \times U \to \Re^+$, $V_2 : \Re^+ \times \mathcal{X} \times U \to \Re^+$, *with* $V_i(t,0,0) = 0$ *for all* $t \geq 0$ ( $i = 1,2$ ) *such that the following hypotheses hold:*

**(H1)** *There exist functions* $\sigma_1 \in KL$, $\beta_1, \mu_1, c_1, \delta_1, \delta_1^u, q_1^u \in K^+$, $\gamma_1, \gamma_1^u, a_1, p_1, p_1^u \in N$, $L_1 : \Re^+ \times \mathcal{X} \to \Re^+$ *with* $L_1(t,0) = 0$ *for all* $t \geq 0$, *such that for every* $(t_0, x_0, u, d) \in \Re^+ \times \mathcal{X} \times M_U \times M_D$ *the mapping* $t \to V_1(t, \phi(t, t_0, x_0, u, d), u(t))$ *is locally bounded on* $[t_0, t_{\max})$ *and the following estimates hold for all* $t \in [t_0, t_{\max})$:

$$V_1(t, \phi(t, t_0, x_0, u, d), u(t)) \leq \sigma_1\big(\beta_1(t_0) L_1(t_0, x_0), t - t_0\big) + \sup_{t_0 \leq \tau \leq t} \gamma_1\big(\delta_1(\tau) V_2(\tau)\big) + \sup_{t_0 \leq \tau \leq t} \gamma_1^u\big(\delta_1^u(\tau) \|u(\tau)\|_{\mathcal{U}}\big) \quad (3.1)$$

$$\beta_1(t) L_1(t, \phi(t, t_0, x_0, u, d)) \leq \max\left\{ \mu_1(t - t_0), c_1(t_0), a_1(\|x_0\|_{\mathcal{X}}), \sup_{t_0 \leq \tau \leq t} p_1(V_2(\tau)), \sup_{t_0 \leq \tau \leq t} p_1^u\big(q_1^u(\tau) \|u(\tau)\|_{\mathcal{U}}\big) \right\} \quad (3.2)$$

*where* $V_2(t) = V_2\big(t, \phi(t, t_0, x_0, u, d), u(t)\big)$ *and* $t_{\max}$ *is the maximal existence time of the transition map of* $\Sigma$.

**(H2)** *There exist functions* $\sigma_2 \in KL$, $\beta_2, \mu_2, c_2, \delta_2, \delta_2^u, q_2^u \in K^+$, $\gamma_2, \gamma_2^u, a_2, p_2, p_2^u \in N$, $L_2 : \Re^+ \times \mathcal{X} \to \Re^+$ *with* $L_2(t,0) = 0$ *for all* $t \geq 0$, *such that for every* $(t_0, x_0, u, d) \in \Re^+ \times \mathcal{X} \times M_U \times M_D$ *the mapping* $t \to V_2(t, \phi(t, t_0, x_0, u, d), u(t))$ *is locally bounded on* $[t_0, t_{\max})$ *and the following estimates hold for all* $t \in [t_0, t_{\max})$:

$$V_2(t, \phi(t, t_0, x_0, u, d), u(t)) \leq \sigma_2\big(\beta_2(t_0) L_2(t_0, x_0), t - t_0\big) + \sup_{t_0 \leq \tau \leq t} \gamma_2\big(\delta_2(\tau) V_1(\tau)\big) + \sup_{t_0 \leq \tau \leq t} \gamma_2^u\big(\delta_2^u(\tau) \|u(\tau)\|_{\mathcal{U}}\big) \quad (3.3)$$

$$\beta_2(t) L_2(t, \phi(t, t_0, x_0, u, d)) \leq \max\left\{ \mu_2(t - t_0), c_2(t_0), a_2(\|x_0\|_{\mathcal{X}}), \sup_{t_0 \leq \tau \leq t} p_2(V_1(\tau)), \sup_{t_0 \leq \tau \leq t} p_2^u\big(q_2^u(\tau) \|u(\tau)\|_{\mathcal{U}}\big) \right\} \quad (3.4)$$

*where* $V_1(t) = V_1\big(t, \phi(t, t_0, x_0, u, d), u(t)\big)$ *and* $t_{\max}$ *is the maximal existence time of the transition map of* $\Sigma$.

**(H3)** *There exist a function* $\rho \in K_\infty$ *and a constant* $M > 0$ *such that:*

$$\delta_1(t) \leq M, \quad \forall t \geq 0 \quad (3.5)$$

$$g_1\big(\delta_1(t) g_2(\delta_2(\tau) s)\big) \leq s, \quad \forall t, s \geq 0 \text{ and } \tau \in [0, t] \quad (3.6)$$

*where* $g_i(s) := \gamma_i(s) + \rho(\gamma_i(s))$, $i = 1,2$.

**(H4)** *There exists a function* $a \in N$ *such that the following inequality holds for all* $(t, x, u) \in \Re^+ \times \mathcal{X} \times U$:

$$\|H(t, x, u)\|_Y \leq a\big(V_1(t, x, u) + \gamma_1(\delta_1(t) V_2(t, x, u))\big) \quad (3.7)$$

**(H5)** *There exists a function* $b \in N$ *such that the following inequality holds for all* $(t, x) \in \Re^+ \times \mathcal{X}$:

$$\|x\|_{\mathcal{X}} \leq b\big(L_1(t, x) + L_2(t, x)\big); \quad \max(L_1(t, x), L_2(t, x)) \leq b(\|x\|_{\mathcal{X}}) \quad (3.8)$$

*Then there exists a function* $\gamma \in N$ *such that system* $\Sigma$ *satisfies the WIOS property from the input* $u \in M_U$ *with gain* $\gamma \in N$ *and weight* $\delta \in K^+$, *where*

$$\delta(t) := \max\{\delta_1^u(t), \delta_2^u(t), q_1^u(t), q_2^u(t)\} \quad (3.9)$$

*Moreover, if* $\beta_1, \beta_2, c_1, c_2, \delta_2 \in K^+$ *are bounded then system* $\Sigma$ *satisfies the UWIOS property from the input* $u \in M_U$ *with gain* $\gamma \in N$ *and weight* $\delta \in K^+$.

**Remark 3.2:**

(a) It should be clear that Theorem 3.1 takes into account all possible cases (weights, non-uniformity with respect to initial times) and thus is applicable to a very wide class of systems.

(b) When $\gamma_1 \in N$ (or $\gamma_2 \in N$) is identically zero, it follows that (3.6) is automatically satisfied. This is the case of systems in cascade (see [10]). On the other hand, if $\gamma_i(s) = K_i s$ for certain constants $K_i \geq 0$ ($i = 1,2$) then inequality (3.6) is satisfied if $K_1 K_2 \sup_{t \geq 0} \left( \delta_1(t) \max_{\tau \in [0,t]} \delta_2(\tau) \right) < 1$. Moreover, if $\gamma_i(s) = K_i s$ for certain constants $K_i \geq 0$ ($i = 1,2$) and $\delta_1(t) \equiv \delta_2(t) \equiv 1$ then hypothesis (H3) is satisfied if $K_1 K_2 < 1$. This is the case of the classical Small-Gain Theorem.

(c) If, instead of hypothesis (H4), there exists a function $a \in N$ such that $\|H(t,x,u)\|_Y \leq a\left( V_2(t,x,u) + \gamma_2(\delta_2(t) V_1(t,x,u)) \right)$ holds for all $(t,x,u) \in \Re^+ \times X \times U$, then indices 1 and 2 must be changed in hypothesis (H3). Furthermore, in this case system $\Sigma$ satisfies the UWIOS property from the input $u \in M_U$ with gain $\gamma \in N$ and weight $\delta \in K^+$, if functions $\beta_1, \beta_2, c_1, c_2, \delta_1 \in K^+$ are bounded.

(d) In classical Small-Gain Theorems (e.g. [10]), the functions $L_1 : \Re^+ \times X \to \Re^+$ and $L_2 : \Re^+ \times X \to \Re^+$ take the form $L_1(t,x) = x_1$ and $L_2(t,x) = x_2$, where $x = (x_1, x_2)$. It follows that hypothesis (H5) automatically holds with $b(s) := s$. This is the case in Corollary 3.4 below.

Since Small-Gain results are frequently applied to feedback interconnections of control systems, we need to clarify the notion of the feedback interconnection of two control systems. However, the fact that we do not require the classical semigroup property for each of the interconnected subsystems, creates technical difficulties: for example, the determination of the set of sampling times for the composite system is not trivial. In order to guarantee the existence of a set of sampling times for the composite system, we assume that the sampling times of the composite system are the common sampling times of the interconnected subsystems. The details are given in the following definition.

**Definition 3.3:** *Consider a pair of control systems $\Sigma_1 = (X_1, Y_1, M_{S_2 \times U}, M_D, \tilde{\phi}_1, \pi_1, H_1)$, $\Sigma_2 = (X_2, Y_2, M_{S_1 \times U}, M_D, \tilde{\phi}_2, \pi_2, H_2)$ with outputs $H_1 : \Re^+ \times X_1 \times U \to S_1 \subseteq Y_1$, $H_2 : \Re^+ \times X_2 \times Y_1 \times U \to S_2 \subseteq Y_2$ and the BIC property and for which $0 \in X_i$ $i = 1,2$ are robust equilibrium points from the inputs $(v_2, u) \in M_{S_2 \times U}$, $(v_1, u) \in M_{S_1 \times U}$, respectively. Suppose that there exists a **unique** pair of a map $\phi = (\phi_1, \phi_2) : A_\phi \to X$ and a set-valued map $\Re^+ \times X \times M_U \times M_D \ni (t_0, x_0, u, d) \to \pi(t_0, x_0, u, d) \subseteq [t_0, +\infty)$, where $A_\phi \subseteq \Re^+ \times \Re^+ \times X \times M_U \times M_D$, $X = X_1 \times X_2$, such that for every $(t, t_0, x_0, u, d) \in A_\phi$ with $t \geq t_0$, $x_0 = (x_{1,0}, x_{2,0}) \in X_1 \times X_2$ it holds that:*

"*there exists a pair of external inputs $v_i \subseteq M(S_i)$ $i = 1,2$ with $v_1(\tau) = H_1(\tau, \phi_1(\tau, t_0, x_0, u, d), u(\tau))$, $v_2(\tau) = H_2(\tau, \phi_2(\tau, t_0, x_0, u, d), v_1(\tau), u(\tau))$ for all $\tau \in [t_0, t]$, $(v_i, u) \in M_{S_i \times U}$ $i = 1,2$, $\pi(t_0, x_0, u, d) = \pi_1(t_0, x_{1,0}, (v_2, u), d) \cap \pi_2(t_0, x_{2,0}, (v_1, u), d)$ and $\phi_1(\tau, t_0, x_0, u, d) = \tilde{\phi}_1(\tau, t_0, x_{1,0}, (v_2, u), d)$, $\phi_2(\tau, t_0, x_0, u, d) = \tilde{\phi}_2(\tau, t_0, x_{2,0}, (v_1, u), d)$ for all $\tau \in [t_0, t]$.*"

*Moreover, let $Y$ be a normed linear space and $H : \Re^+ \times X \times U \to Y$ a continuous map that maps bounded sets of $\Re^+ \times X \times U$ into bounded sets of $Y$, with $H(t,0,0) = 0$ for all $t \geq 0$ and suppose that $\Sigma := (X, Y, M_U, M_D, \phi, \pi, H)$ is a control system with outputs and the BIC property, for which $0 \in X$ is a robust equilibrium point from the input $u \in M_U$. Then system $\Sigma$ is said to be the **feedback connection** or the **interconnection** of systems $\Sigma_1$ and $\Sigma_2$.*

It should be emphasized that the feedback interconnection of two systems may create a system, which has different qualitative properties from each of the interconnected subsystems. For example, if we interconnect a subsystem described by RFDEs (see Example 2.10) with a hybrid subsystem with impulses at fixed times (see Example 2.12), then the overall system will be a system with both "memory" and impulses (discontinuous systems described by RFDEs-see [42]).

We are now in a position to state our main result for feedback interconnections of control systems. It is a direct consequence of Theorem 3.1 and its proof is omitted.

**Corollary 3.4:** *Suppose that $\Sigma := (X, Y, M_U, M_D, \phi, \pi, H)$ is the feedback connection of systems $\Sigma_1 = (X_1, Y_1, M_{S_2 \times U}, M_D, \tilde{\phi}_1, \pi_1, H_1)$ and $\Sigma_2 = (X_2, Y_2, M_{S_1 \times U}, M_D, \tilde{\phi}_2, \pi_2, H_2)$ with outputs $H_1 : \Re^+ \times X_1 \times U \to S_1 \subseteq Y_1$, $H_2 : \Re^+ \times X_2 \times Y_1 \times U \to S_2 \subseteq Y_2$. We assume that:*

**(H1')** *Subsystem $\Sigma_1$ satisfies the WIOS property from the inputs $v_2 \in M(S_2)$ and $u \in M_U$. Particularly, there exist functions $\sigma_1 \in KL$, $\beta_1, \mu_1, c_1, \delta_1, \delta_1^u, q_1^u \in K^+$, $\gamma_1, \gamma_1^u, a_1, p_1, p_1^u \in N$, such that the following estimate holds for all $(t_0, x_1, (v_2, u_0), d) \in \Re^+ \times X_1 \times M_{S_2 \times U} \times M_D$ and $t \geq t_0$:*

$$\left\| H_1(t, \tilde{\phi}_1(t, t_0, x_1, (v_2, u), d), u(t)) \right\|_{Y_1} \leq \sigma_1 \left( \beta_1(t_0) \|x_1\|_{X_1}, t - t_0 \right) + \sup_{t_0 \leq \tau \leq t} \gamma_1 \left( \delta_1(\tau) \|v_2(\tau)\|_{Y_2} \right) + \sup_{t_0 \leq \tau \leq t} \gamma_1^u \left( \delta_1^u(\tau) \|u(\tau)\|_U \right)$$
(3.10)

$$\beta_1(t) \left\| \tilde{\phi}_1(t, t_0, x_1, (v_2, u), d) \right\|_{X_1} \leq \max \left\{ \mu_1(t - t_0), c_1(t_0), a_1 \left( \|x_1\|_{X_1} \right), \sup_{t_0 \leq \tau \leq t} p_1 \left( \|v_2(\tau)\|_{Y_2} \right), \sup_{t_0 \leq \tau \leq t} p_1^u \left( q_1^u(\tau) \|u(\tau)\|_U \right) \right\}$$
(3.11)

**(H2')** *Subsystem $\Sigma_2$ satisfies the WIOS property from the inputs $v_1 \in M(S_1)$ and $u \in M_U$. Particularly, there exist functions $\sigma_2 \in KL$, $\beta_2, \mu_2, c_2, \delta_2, \delta_2^u, q_2^u \in K^+$, $\gamma_2, \gamma_2^u, a_2, p_2, p_2^u \in N$, such that the following estimate holds for all $(t_0, x_2, (v_1, u_0), d) \in \Re^+ \times X_2 \times M_{S_1 \times U} \times M_D$ and $t \geq t_0$:*

$$\left\| H_2(t, \tilde{\phi}_2(t, t_0, x_2, (v_1, u), d), v_1(t), u(t)) \right\|_{Y_2} \leq \sigma_2 \left( \beta_2(t_0) \|x_2\|_{X_2}, t - t_0 \right) + \sup_{t_0 \leq \tau \leq t} \gamma_2 \left( \delta_2(\tau) \|v_1(\tau)\|_{Y_1} \right) + \sup_{t_0 \leq \tau \leq t} \gamma_2^u \left( \delta_2^u(\tau) \|u(\tau)\|_U \right) \quad (3.12)$$

$$\beta_2(t) \left\| \tilde{\phi}_2(t, t_0, x_2, (v_1, u), d) \right\|_{X_2} \leq \max \left\{ \mu_2(t - t_0), c_2(t_0), a_2 \left( \|x_2\|_{X_2} \right), \sup_{t_0 \leq \tau \leq t} p_2 \left( \|v_1(\tau)\|_{Y_1} \right), \sup_{t_0 \leq \tau \leq t} p_2^u \left( q_2^u(\tau) \|u(\tau)\|_U \right) \right\} \quad (3.13)$$

*Moreover, assume that hypothesis (H3) of Theorem 3.1 holds and there exists a function $a \in N$ such that the following inequality holds for all $(t, x, u) \in \Re^+ \times X \times U$ with $x = (x_1, x_2) \in X_1 \times X_2$:*

$$\| H(t, x, u) \|_Y \leq a \left( \| H_1(t, x_1, u) \|_{Y_1} + \gamma_1 \left( \delta_1(t) \| H_2(t, x_2, H_1(t, x_1, u), u) \|_{Y_2} \right) \right)$$
(3.14)

*Then there exists a function $\gamma \in N$ such that system $\Sigma$ satisfies the WIOS property from the input $u \in M_U$ with gain $\gamma \in N$ and weight $\delta \in K^+$, where $\delta \in K^+$ is defined by (3.9). Moreover, if $\beta_1, \beta_2, c_1, c_2, \delta_2 \in K^+$ are bounded then system $\Sigma$ satisfies the UWIOS property from the input $u \in M_U$ with gain $\gamma \in N$ and weight $\delta \in K^+$.*

**Remark 3.5:**

(a) When $\delta_1(t) \equiv \delta_2(t) \equiv 1$ and $\gamma_1 \in K_\infty$ then the result of Corollary 3.4 guarantees the WIOS property from the input $u \in M_U$ for the output $H(t, x, u) := (H_1(t, x_1, u), H_2(t, x_2, H_1(t, x_1, u), u))$, i.e., for the output that combines the outputs of each individual subsystem. Moreover, if in addition the functions $q_i^u(t)$ ($i = 1, 2$) are bounded, then the result of Corollary 3.4 guarantees the IOS from the input $u \in M_U$ for the output $H(t, x, u) := (H_1(t, x_1, u), H_2(t, x_2, H_1(t, x_1, u), u))$. Finally, if in addition the functions $\beta_1, \beta_2, c_1, c_2, \delta_2 \in K^+$ are bounded then the result of Corollary 3.4 guarantees the UIOS from the input $u \in M_U$ for the output $H(t, x, u) := (H_1(t, x_1, u), H_2(t, x_2, H_1(t, x_1, u), u))$. This particular case coincides with the prior result of the nonlinear ISS Small-Gain Theorem presented in [10] for control systems described by ODEs.

**(b)** Conditions (3.11) and (3.13) hold automatically, when each one of the subsystems $\Sigma_1$ and $\Sigma_2$ satisfy the WISS property.

**(c)** If, instead of hypothesis (3.14), there exists a function $a \in N$ such that $\|H(t,x,u)\|_Y \le a\left(\|H_2(t,x_2,H_1(t,x_1,u),u)\|_{Y_2} + \gamma_2\left(\delta_2(t)\|H_1(t,x_1,u)\|_{Y_1}\right)\right)$ holds for all $(t,x,u) \in \Re^+ \times X \times U$, then indices 1 and 2 must be changed in hypothesis (H3). Furthermore, in this case system $\Sigma$ satisfies the UWIOS property from the input $u \in M_U$ with gain $\gamma \in N$ and weight $\delta \in K^+$, if functions $\beta_1, \beta_2, c_1, c_2, \delta_1 \in K^+$ are bounded.

**Proof of Theorem 3.1:** The proof consists of three steps:

<u>Step 1:</u> We show that $\Sigma$ is RFC from the input $u \in M_U$.

<u>Step 2:</u> Let $\tilde{\varphi}(s) := s + \frac{1}{2}\rho(s)$, where $\rho \in K_\infty$ is the function involved in hypothesis (H3). We show that properties P1 and P2 of Lemma 2.16 hold for system $\Sigma$ with $V = V_1$ or $V = \tilde{\varphi}(\gamma_1(\delta_1(t)V_2))$, for appropriate $\tilde{\gamma} \in N$ and $\delta \in K^+$ as defined by (3.9). Moreover, if $\beta_1, \beta_2 \in K^+$ are bounded we show that properties P1 and P2 of Lemma 2.17 hold for system $\Sigma$ with $V = V_1$ or $V = \tilde{\varphi}(\gamma_1(\delta_1(t)V_2))$, for appropriate $\tilde{\gamma} \in N$ and $\delta \in K^+$ as defined by (3.9).

<u>Step 3:</u> We show that property P3 of Lemma 2.16 holds for system $\Sigma$ with $V = V_1$ or $V = \tilde{\varphi}(\gamma_1(\delta_1(t)V_2))$, for appropriate $\tilde{\gamma} \in N$ and $\delta \in K^+$ as defined by (3.9). Moreover, if $\beta_1, \beta_2, c_1, c_2, \delta_2 \in K^+$ are bounded, we show that property P3 of Lemma 2.17 holds for system $\Sigma$ with $V = V_1$ or $V = \tilde{\varphi}(\gamma_1(\delta_1(t)V_2))$, for appropriate $\tilde{\gamma} \in N$ and $\delta \in K^+$ as defined by (3.9).

It then follows from Lemma 2.16 that exist functions $\sigma \in KL$ and $\beta \in K^+$ such that the following estimate holds for all $u \in M_U$, $(t_0, x_0, d) \in \Re^+ \times X \times M_D$ and $t \ge t_0$:

$$V_1(t, \phi(t,t_0,x_0,u,d), u(t)) \le \sigma\left(\beta(t_0)\|x_0\|_X, t-t_0\right) + \sup_{t_0 \le \tau \le t} \gamma\left(\delta(\tau)\|u(\tau)\|_U\right) \quad (3.15a)$$

$$\tilde{\varphi}\left(\gamma_1(\delta_1(t)V_2(t,\phi(t,t_0,x_0,u,d),u(t)))\right) \le \sigma\left(\beta(t_0)\|x_0\|_X, t-t_0\right) + \sup_{t_0 \le \tau \le t} \gamma\left(\delta(\tau)\|u(\tau)\|_U\right) \quad (3.15b)$$

Moreover, if $\beta_1, \beta_2, c_1, c_2, \delta_2 \in K^+$ are bounded it follows from Lemma 2.17 that estimates (3.15a,b) hold with $\beta(t) \equiv 1$.

Thus, using (3.7) and the fact that $\tilde{\varphi}(s) \ge s$ for all $s \ge 0$, we conclude that $\Sigma$ satisfies the WIOS property from the input $u \in M_U$ with gain $\gamma(s) := a(4\tilde{\gamma}(s)) \in N$ and weight $\delta \in K^+$. Moreover, if $\beta_1, \beta_2, c_1, c_2, \delta_2 \in K^+$ are bounded we conclude that $\Sigma$ satisfies the UWIOS property from the input $u \in M_U$ with gain $\gamma(s) := a(4\tilde{\gamma}(s)) \in N$ and weight $\delta \in K^+$.

<u>Step 1:</u>

Let arbitrary $(t,t_0,x_0,u,d) \in A_\phi$ with $t \ge t_0$, $x_0 \in X$ and let $t_{\max} \in (t_0, +\infty]$ the maximal existence time of the transition map $\phi$ of $\Sigma := (X, Y, M_U, M_D, \phi, \pi, H)$ that corresponds to $(t_0,x_0,u,d) \in \Re^+ \times X \times M_U \times M_D$. Notice that by virtue of the BIC property, if $t_{\max} < +\infty$ then for every $M > 0$ there exists $t \in [t_0, t_{\max})$ with $\|\phi(t,t_0,x_0,u,d)\|_\mathcal{X} > M$. We define $V_1(\tau) = V_1(\tau,\phi(\tau,t_0,x_0,u,d),u(\tau))$, $L_1(\tau) = L_1(\tau,\phi(\tau,t_0,x_0,u,d))$ and $V_2(\tau) = V_2(\tau,\phi(\tau,t_0,x_0,u,d),u(\tau))$, $L_2(\tau) = L_2(\tau,\phi(\tau,t_0,x_0,u,d))$ for all $\tau \in [t_0,t]$.

The previous definitions in conjunction with (3.1), (3.2), (3.3), (3.4) imply the following inequalities for all $t \in [t_0,t_{\max})$:

$$V_1(t) \le \sigma_1\big(\beta_1(t_0)L_1(t_0), t-t_0\big) + \sup_{t_0 \le \tau \le t} \gamma_1\big(\delta_1(\tau)V_2(\tau)\big) + \sup_{t_0 \le \tau \le t} \gamma_1^u\big(\delta_1^u(\tau)\|u(\tau)\|_\mathcal{U}\big) \qquad (3.16)$$

$$\beta_1(t)L_1(t) \le \max\left\{\mu_1(t-t_0), c_1(t_0), a_1(\|x_0\|_\mathcal{X}), \sup_{t_0 \le \tau \le t} p_1(V_2(\tau)), \sup_{t_0 \le \tau \le t} p_1^u\big(q_1^u(\tau)\|u(\tau)\|_\mathcal{U}\big)\right\} \qquad (3.17)$$

$$V_2(t) \le \sigma_2\big(\beta_2(t_0)L_2(t_0), t-t_0\big) + \sup_{t_0 \le \tau \le t} \gamma_2\big(\delta_2(\tau)V_1(\tau)\big) + \sup_{t_0 \le \tau \le t} \gamma_2^u\big(\delta_2^u(\tau)\|u(\tau)\|_\mathcal{U}\big) \qquad (3.18)$$

$$\beta_2(t)L_2(t) \le \max\left\{\mu_2(t-t_0), c_2(t_0), a_2(\|x_0\|_\mathcal{X}), \sup_{t_0 \le \tau \le t} p_2(V_1(\tau)), \sup_{t_0 \le \tau \le t} p_2^u\big(q_2^u(\tau)\|u(\tau)\|_\mathcal{U}\big)\right\} \qquad (3.19)$$

Let $\rho \in K_\infty$ the function involved in hypothesis (H3) and define $\kappa(s) := s + \rho^{-1}(s)$, $\varphi(s) := s + \rho(s)$. Using the inequality $r+s \le \max\{\kappa(r); \varphi(s)\}$ (which holds for all $r,s \ge 0$) as well as the equality $g_2(s) = \varphi(\gamma_2(s))$, we obtain from (3.18) for all $t \in [t_0, t_{\max})$:

$$V_2(t) \le \max\left\{\kappa\left(\sigma_2\big(\beta_2(t_0)L_2(t_0), t-t_0\big) + \sup_{t_0 \le \tau \le t} \gamma_2^u\big(\delta_2^u(\tau)\|u(\tau)\|_\mathcal{U}\big)\right); \sup_{t_0 \le \tau \le t} g_2\big(\delta_2(\tau)V_1(\tau)\big)\right\} \qquad (3.20)$$

Notice that inequality (3.6) implies that $\gamma_1\big(\delta_1(t)g_2(\delta_2(\tau)s)\big) \le \varphi^{-1}(s)$, $\forall t, s \ge 0$ and $\tau \in [0,t]$. Thus (3.20) in conjunction with (3.5) and the previous observation implies the following estimate which holds for all $t \in [t_0, t_{\max})$:

$$\gamma_1\big(\delta_1(t)V_2(t)\big) \le \max\left\{\gamma_1\left(M\kappa\left(\sigma_2\big(\beta_2(t_0)L_2(t_0), t-t_0\big) + \sup_{t_0 \le \tau \le t} \gamma_2^u\big(\delta_2^u(\tau)\|u(\tau)\|_\mathcal{U}\big)\right)\right); \varphi^{-1}\left(\sup_{t_0 \le \tau \le t} V_1(\tau)\right)\right\} \qquad (3.21)$$

Combining estimate (3.16) with (3.21), we obtain:

$$\sup_{t_0 \le \tau \le t} V_1(\tau) \le \sigma_1\big(\beta_1(t_0)L_1(t_0), 0\big) + \sup_{t_0 \le \tau \le t} \gamma_1^u\big(\delta_1^u(\tau)\|u(\tau)\|_\mathcal{U}\big)$$
$$+ \max\left\{\gamma_1\left(M\kappa\left(\sigma_2\big(\beta_2(t_0)L_2(t_0), 0\big) + \sup_{t_0 \le \tau \le t} \gamma_2^u\big(\delta_2^u(\tau)\|u(\tau)\|_\mathcal{U}\big)\right)\right); \varphi^{-1}\left(\sup_{t_0 \le \tau \le t} V_1(\tau)\right)\right\} \qquad (3.22)$$

Distinguishing the cases $\gamma_1\left(M\kappa\left(\sigma_2\big(\beta_2(t_0)L_2(t_0), 0\big) + \sup_{t_0 \le \tau \le t} \gamma_2^u\big(\delta_2^u(\tau)\|u(\tau)\|_\mathcal{U}\big)\right)\right) \ge \varphi^{-1}\left(\sup_{t_0 \le \tau \le t} V_1(\tau)\right)$, $\gamma_1\left(M\kappa\left(\sigma_2\big(\beta_2(t_0)L_2(t_0), 0\big) + \sup_{t_0 \le \tau \le t} \gamma_2^u\big(\delta_2^u(\tau)\|u(\tau)\|_\mathcal{U}\big)\right)\right) \le \varphi^{-1}\left(\sup_{t_0 \le \tau \le t} V_1(\tau)\right)$, using the identity $s - \varphi^{-1}(s) = \kappa^{-1}(s)$ and the fact that $\kappa(s) \ge s$ in conjunction with (3.22) and (3.8) (which implies $L_i(t_0) \le b(\|x_0\|_\mathcal{X})$, $i = 1, 2$) gives the following estimate which holds for all $t \in [t_0, t_{\max})$:

$$\sup_{t_0 \leq \tau \leq t} V_1(\tau) \leq$$

$$\max \left\{ \begin{array}{l} \kappa\left(\sigma_1\left(\beta_1(t_0) b\left(\|x_0\|_X\right), 0\right) + \sup_{t_0 \leq \tau \leq t} \gamma_1^u\left(\delta_1^u(\tau) \|u(\tau)\|_U\right)\right) \\ \sigma_1\left(\beta_1(t_0) b\left(\|x_0\|_X\right), 0\right) + \sup_{t_0 \leq \tau \leq t} \gamma_1^u\left(\delta_1^u(\tau) \|u(\tau)\|_U\right) + \gamma_1\left(M\kappa\left(\sigma_2\left(\beta_2(t_0) b\left(\|x_0\|_X\right), 0\right) + \sup_{t_0 \leq \tau \leq t} \gamma_2^u\left(\delta_2^u(\tau) \|u(\tau)\|_U\right)\right)\right) \end{array} \right\}$$

(3.23)

We show next that $\Sigma$ is RFC from the input $u \in M_U$ by contradiction. Suppose that $t_{\max} < +\infty$. Then by virtue of the BIC property for every $M > 0$ there exists $t \in [t_0, t_{\max})$ with $\|\phi(t, t_0, x_0, u, d)\|_X > M$. On the other hand estimate (3.23) in conjunction with the hypothesis $t_{\max} < +\infty$ shows that there exists $M_1 \geq 0$ such that $\sup_{t_0 \leq \tau < t_{\max}} V_1(\tau) \leq M_1$. The fact that $V_1(t)$ is bounded in conjunction with estimates (3.18) and (3.19), implies that there exist constants $M_2, M_3 \geq 0$ such that $\sup_{t_0 \leq \tau < t_{\max}} V_2(\tau) \leq M_2$ and $\sup_{t_0 \leq \tau < t_{\max}} L_2(\tau) \leq M_3$. Finally, the fact that $V_2(t)$ is bounded in conjunction with estimate (3.17), implies that there exists a constant $M_4 \geq 0$ such that $\sup_{t_0 \leq \tau < t_{\max}} L_1(\tau) \leq M_4$. It follows from (3.8) and inequality $\|\phi(t, t_0, x_0, u, d)\|_X \leq b(L_1(t) + L_2(t))$ that the transition map of $\Sigma$, i.e., $\phi(t, t_0, x_0, u, d)$, is bounded on $[t_0, t_{\max})$ and this contradicts the requirement that for every $M > 0$ there exists $t \in [t_0, t_{\max})$ with $\|\phi(t, t_0, x_0, u, d)\|_X > M$. Hence, we must have $t_{\max} = +\infty$.

Let arbitrary $R \geq 0$, $T \geq 0$. For every $u \in \mathcal{M}(B_U[0, R]) \cap M_U$, $s \in [0, T]$, $\|x_0\|_X \leq R$, $t_0 \in [0, T]$, $d \in M_D$ estimate (3.23) shows that there exists $M_1(T, R) \geq 0$ such that $V_1(t_0 + s) \leq M_1(T, R)$, for all $s \in [0, T]$. The previous observation in conjunction with estimates (3.18), (3.19) and (3.8) (which gives $L_i(t_0) \leq b(\|x_0\|_X)$, $i = 1, 2$), implies that there exist $M_2(T, R), M_3(T, R) \geq 0$ such that for every $u \in \mathcal{M}(B_U[0, R]) \cap M_U$, $s \in [0, T]$, $\|x_0\|_X \leq R$, $t_0 \in [0, T]$, $d \in M_D$ we have $V_2(t_0 + s) \leq M_2(T, R)$ and $L_2(t_0 + s) \leq M_3(T, R)$, for all $s \in [0, T]$. Finally, inequality $V_2(t_0 + s) \leq M_2(T, R)$ in conjunction with estimate (3.17), implies that there exists a constant $M_4(T, R) \geq 0$ such that for every $u \in \mathcal{M}(B_U[0, R]) \cap M_U$, $s \in [0, T]$, $\|x_0\|_X \leq R$, $t_0 \in [0, T]$, $d \in M_D$ we have $L_1(t_0 + s) \leq M_4(T, R)$, for all $s \in [0, T]$. It follows from (3.8) and inequality $\|\phi(t, t_0, x_0, u, d)\|_X \leq b(L_1(t) + L_2(t))$ that for every $u \in \mathcal{M}(B_U[0, R]) \cap M_U$, $s \in [0, T]$, $\|x_0\|_X \leq R$, $t_0 \in [0, T]$, $d \in M_D$ the transition map of $\Sigma$, i.e., $\phi(t, t_0, x_0, u, d)$, satisfies $\|\phi(t, t_0, x_0, u, d)\|_X \leq b(M_3(T, R) + M_4(T, R)) < +\infty$ and this according to Definition 2.2 implies that $\Sigma$ is RFC from the input $u \in M_U$.

Step 2:

Using (3.23) in conjunction with the inequality $q(r + s) \leq q(\kappa(r)) + q(\varphi(s))$ (which holds for all $r, s \geq 0$ and $q \in N$) gives the following estimate, which holds for all $t \geq t_0$:

$$\begin{aligned} V_1(t) &\leq \kappa\left(\kappa\left(\sigma_1\left(\beta_1(t_0) b\left(\|x_0\|_X\right), 0\right)\right)\right) + \gamma_1\left(M\kappa\left(\kappa\left(\sigma_2\left(\beta_2(t_0) b\left(\|x_0\|_X\right), 0\right)\right)\right)\right) \\ &+ \sup_{t_0 \leq \tau \leq t} \kappa\left(\varphi\left(\gamma_1^u\left(\delta_1^u(\tau) \|u(\tau)\|_U\right)\right)\right) + \sup_{t_0 \leq \tau \leq t} \gamma_1\left(M\kappa\left(\varphi\left(\gamma_2^u\left(\delta_2^u(\tau) \|u(\tau)\|_U\right)\right)\right)\right) \end{aligned}$$

(3.24)

Moreover, combining estimates (3.21) and (3.23) and using the equalities $\varphi^{-1}(\kappa(s)) = \rho(s)$ and $\varphi(s) := s + \rho(s)$ as well as the inequalities $\varphi^{-1}(s) \leq s$ and (3.8) (which gives $L_i(t_0) \leq b(\|x_0\|_X)$, $i = 1, 2$), gives the following estimate, which holds for all $t \geq t_0$:

$$\gamma_1(\delta_1(t)V_2(t)) \leq \gamma_1\left(M\kappa\left(\sigma_2(\beta_2(t_0)b(\|x_0\|_\mathcal{X}),0) + \sup_{t_0 \leq \tau \leq t}\gamma_2^u(\delta_2^u(\tau)\|u(\tau)\|_\mathcal{U})\right)\right)$$
$$+ \varphi\left(\sigma_1(\beta_1(t_0)b(\|x_0\|_\mathcal{X}),0) + \sup_{t_0 \leq \tau \leq t}\gamma_1^u(\delta_1^u(\tau)\|u(\tau)\|_\mathcal{U})\right) \quad (3.25)$$

Using (3.25) in conjunction with the inequality $q(r+s) \leq q(\kappa(r)) + q(\varphi(s))$ (which holds for all $r,s \geq 0$ and $q \in N$) gives the following estimate, which holds for all $t \geq t_0$:

$$\gamma_1(\delta_1(t)V_2(t)) \leq \gamma_1\left(M\kappa\left(\kappa(\sigma_2(\beta_2(t_0)b(\|x_0\|_\mathcal{X}),0))\right)\right) + \sup_{t_0 \leq \tau \leq t}\gamma_1\left(M\kappa\left(\varphi(\gamma_2^u(\delta_2^u(\tau)\|u(\tau)\|_\mathcal{U}))\right)\right)$$
$$+ \varphi\left(\kappa(\sigma_1(\beta_1(t_0)b(\|x_0\|_\mathcal{X}),0))\right) + \sup_{t_0 \leq \tau \leq t}\varphi\left(\varphi(\gamma_1^u(\delta_1^u(\tau)\|u(\tau)\|_\mathcal{U}))\right) \quad (3.26)$$

Let $\tilde{\varphi}(s) := s + \frac{1}{2}\rho(s)$. Using (3.26) in conjunction with the inequality $q(r+s) \leq q(b(r)) + q(\varphi(s))$ (which holds for all $r,s \geq 0$ and $q \in N$), we obtain the following estimate, which holds for all $t \geq t_0$:

$$\tilde{\varphi}(\gamma_1(\delta_1(t)V_2(t))) \leq \tilde{\varphi}\left(\kappa\left(\gamma_1\left(M\kappa\left(\kappa(\sigma_2(\beta_2(t_0)b(\|x_0\|_\mathcal{X}),0))\right)\right) + \varphi(\kappa(\sigma_1(\beta_1(t_0)b(\|x_0\|_\mathcal{X}),0)))\right)\right)$$
$$+ \sup_{t_0 \leq \tau \leq t}\tilde{\varphi}\left(\varphi\left(\varphi(\gamma_1^u(\delta_1^u(\tau)\|u(\tau)\|_\mathcal{U}))\right) + \gamma_1\left(M\kappa\left(\varphi(\gamma_2^u(\delta_2^u(\tau)\|u(\tau)\|_\mathcal{U}))\right)\right)\right) \quad (3.27)$$

Estimates (3.24), (3.27) show that properties P1 and P2 of Lemma 2.16 hold for system $\Sigma$ with $V = V_1$ or $V = \tilde{\varphi}(\gamma_1(\delta_1(t)V_2))$, for appropriate $\tilde{\gamma} \in N$ and $\delta \in K^+$ as defined by (3.9). Moreover, if $\beta_1, \beta_2 \in K^+$ are bounded then estimates (3.24), (3.27) show that properties P1 and P2 of Lemma 2.17 hold for system $\Sigma$ with $V = V_1$ or $V = \tilde{\varphi}(\gamma_1(\delta_1(t)V_2))$, for appropriate $\tilde{\gamma} \in N$ and $\delta \in K^+$ as defined by (3.9). Particularly, $\tilde{\gamma} \in N$ satisfies

$$\tilde{\gamma}(s) \geq \tilde{\varphi}\left(\varphi(\varphi(\gamma_1^u(s))) + \gamma_1(M\kappa(\varphi(\gamma_2^u(s))))\right) \text{ and } \tilde{\gamma}(s) \geq \kappa(\varphi(\gamma_1^u(s))) + \gamma_1(M\kappa(\varphi(\gamma_2^u(s)))), \text{ for all } s \geq 0 \quad (3.28)$$

Step 3:

Let $\tilde{\varphi}(s) := s + \frac{1}{2}\rho(s)$, $\tilde{\kappa}(s) := s + \rho^{-1}(2s)$. Exploiting estimates (3.16), (3.21) in conjunction with the inequality $r + s \leq \max\{\tilde{\kappa}(r); \tilde{\varphi}(s)\}$ (which holds for all $r, s \geq 0$), we obtain:

$$V_1(t) \leq \max\left\{\tilde{\kappa}\left(\sigma_1(\beta_1(t_0)L_1(t_0), t-t_0) + \sup_{t_0 \leq \tau \leq t}\gamma_1^u(\delta_1^u(\tau)\|u(\tau)\|_\mathcal{U})\right); \sup_{t_0 \leq \tau \leq t}\tilde{\varphi}(\gamma_1(\delta_1(\tau)V_2(\tau)))\right\} \quad (3.29)$$

$$\tilde{\varphi}(\gamma_1(\delta_1(t)V_2(t))) \leq \max\left\{\tilde{\varphi}\left(\gamma_1\left(M\kappa\left(\sigma_2(\beta_2(t_0)L_2(t_0), t-t_0) + \sup_{t_0 \leq \tau \leq t}\gamma_2^u(\delta_2^u(\tau)\|u(\tau)\|_\mathcal{U})\right)\right)\right); \tilde{\varphi}\left(\varphi^{-1}\left(\sup_{t_0 \leq \tau \leq t}V_1(\tau)\right)\right)\right\}$$
$$\quad (3.30)$$

Let arbitrary $\xi \in \pi(t_0, x_0, u, d)$ and $t \geq \xi$. Estimates (3.29), (3.30) in conjunction with estimates (3.17), (3.19) and the weak semigroup property imply:

$$V_1(t) \leq \max\left\{\begin{array}{l}\tilde{\kappa}(\kappa(\sigma_1(a_1(\|x_0\|_\mathcal{X}) + c_1(t_0) + \mu_1(\xi-t_0), t-\xi))) \\ \sup_{t_0 \leq \tau \leq \xi}\tilde{\kappa}(\kappa(\sigma_1(p_1^u(q_1^u(\tau)\|u(\tau)\|_\mathcal{U}),0))), \sup_{t_0 \leq \tau \leq \xi}\tilde{\kappa}(\kappa(\sigma_1(p_1(V_2(\tau)), t-\xi))) \\ \sup_{\xi \leq \tau \leq t}\tilde{\varphi}(\gamma_1(\delta_1(\tau)V_2(\tau))), \sup_{\xi \leq \tau \leq t}\tilde{\kappa}(\varphi(\gamma_1^u(\delta_1^u(\tau)\|u(\tau)\|_\mathcal{U})))\end{array}\right\} \quad (3.32)$$

$$\tilde{\varphi}(\gamma_1(\delta_1(t)V_2(t))) \leq \max \begin{cases} \tilde{\varphi}(\gamma_1(M\kappa(\kappa(\sigma_2(a_2(\|x_0\|_X)+c_2(t_0)+\mu_2(\xi-t_0),t-\xi))))) \\ \sup_{\xi\leq\tau\leq t} \tilde{\varphi}(\gamma_1(M\kappa(\varphi(\gamma_2^u(\delta_2^u(\tau)\|u(\tau)\|_U))))), \sup_{t_0\leq\tau\leq\xi} \tilde{\varphi}(\gamma_1(M\kappa(\kappa(\sigma_2(p_2^u(q_2^u(\tau)\|u(\tau)\|_U),0))))) \\ \tilde{\varphi}\left(\varphi^{-1}\left(\sup_{\xi\leq\tau\leq t} V_1(\tau)\right)\right), \sup_{t_0\leq\tau\leq\xi} \tilde{\varphi}(\gamma_1(M\kappa(\kappa(\sigma_2(p_2(V_1(\tau)),t-\xi)))))\end{cases} \quad (3.33)$$

Estimate (3.33) combined with estimate (3.24) gives:

$$\tilde{\varphi}(\gamma_1(\delta_1(t)V_2(t))) \leq \max \begin{cases} \tilde{\varphi}(\gamma_1(M\kappa(\kappa(\sigma_2(a_2(\|x_0\|_X)+c_2(t_0)+\mu_2(\xi-t_0)+p_2(\kappa(A))),t-\xi))))) \\ \sup_{\xi\leq\tau\leq t} \tilde{\varphi}(\gamma_1(M\kappa(\varphi(\gamma_2^u(\delta_2^u(\tau)\|u(\tau)\|_U))))), \sup_{t_0\leq\tau\leq\xi} \tilde{\varphi}(\gamma_1(M\kappa(\kappa(\sigma_2(p_2^u(q_2^u(\tau)\|u(\tau)\|_U),0))))) \\ \tilde{\varphi}\left(\varphi^{-1}\left(\sup_{\xi\leq\tau\leq t} V_1(\tau)\right)\right), \sup_{t_0\leq\tau\leq t} \tilde{\varphi}(\gamma_1(M\kappa(\kappa(\sigma_2(p_2(\varphi(B)),0))))) \end{cases} \quad (3.34)$$

where

$$A = \kappa(\kappa(\sigma_1(\beta_1(t_0)b(\|x_0\|_X),0)))+\gamma_1(M\kappa(\kappa(\sigma_2(\beta_2(t_0)b(\|x_0\|_X),0))))$$
$$B = \kappa(\varphi(\gamma_1^u(\delta_1^u(\tau)\|u(\tau)\|_U)))+\gamma_1(M\kappa(\varphi(\gamma_2^u(\delta_2^u(\tau)\|u(\tau)\|_U))))$$

Similarly, estimate (3.18) combined with estimate (3.24) and inequality (3.8) (which implies $L_i(t_0) \leq b(\|x_0\|_X)$, $i=1,2$), gives:

$$\begin{aligned} V_2(t) &\leq \sigma_2(\beta_2(t_0)b(\|x_0\|_X),0) + \sup_{t_0\leq\tau\leq t}\gamma_2^u(\delta_2^u(\tau)\|u(\tau)\|_U) \\ &+ \gamma_2(\tilde{\delta}_2(t)\,\kappa(\kappa(\kappa(\sigma_1(\beta_1(t_0)b(\|x_0\|_X),0)))+\gamma_1(M\kappa(\kappa(\sigma_2(\beta_2(t_0)b(\|x_0\|_X),0)))))) \\ &+ \sup_{t_0\leq\tau\leq t}\gamma_2(\tilde{\delta}_2(t)\,\varphi(\kappa(\varphi(\gamma_1^u(\delta_1^u(\tau)\|u(\tau)\|_U)))+\gamma_1(M\kappa(\varphi(\gamma_2^u(\delta_2^u(\tau)\|u(\tau)\|_U)))))) \end{aligned} \quad (3.35)$$

where

$$\tilde{\delta}_2(t) := \max_{0\leq\tau\leq t}\delta_2(\tau)$$

Consequently, by combining estimates (3.32) and (3.35) we obtain:

$$V_1(t) \leq \max \begin{cases} \tilde{\kappa}(\kappa(\sigma_1(a_1(\|x_0\|_X)+c_1(t_0)+\mu_1(\xi-t_0),t-\xi))), \tilde{\kappa}(\kappa(\sigma_1(p_1(\kappa(C)),t-\xi))) \\ \sup_{t_0\leq\tau\leq\xi}\tilde{\kappa}(\kappa(\sigma_1(p_1^u(q_1^u(\tau)\|u(\tau)\|_U),0))), \sup_{t_0\leq\tau\leq t}\tilde{\kappa}(\kappa(\sigma_1(p_1(\varphi(D)),0))) \\ \sup_{\xi\leq\tau\leq t}\tilde{\varphi}(\gamma_1(\delta_1(\tau)V_2(\tau))), \sup_{\xi\leq\tau\leq t}\tilde{\kappa}(\varphi(\gamma_1^u(\delta_1^u(\tau)\|u(\tau)\|_U))), \tilde{\kappa}(\kappa(\sigma_1(p_1(E),t-\xi))) \end{cases} \quad (3.36)$$

where

$$C :=$$
$$\sigma_2(\beta_2(t_0)b(\|x_0\|_X),0)+\gamma_2(\tilde{\delta}_2(\xi)\,\kappa(\kappa(\kappa(\sigma_1(\beta_1(t_0)b(\|x_0\|_X),0)))+\gamma_1(M\kappa(\kappa(\sigma_2(\beta_2(t_0)b(\|x_0\|_X),0))))))$$
$$D := \gamma_2^u(\delta_2^u(\tau)\|u(\tau)\|_U)+\gamma_2\left(\varphi\left(\frac{1}{2}\varphi^2(\kappa(\varphi(\gamma_1^u(\delta_1^u(\tau)\|u(\tau)\|_U)))+\gamma_1(M\kappa(\varphi(\gamma_2^u(\delta_2^u(\tau)\|u(\tau)\|_U)))))\right)\right)$$
$$E := \gamma_2\left(\kappa\left(\frac{1}{2}\tilde{\delta}_2^2(\xi)\right)\right)$$

From (3.34) and (3.36) we conclude that there exist functions $S_1, S_2 \in KL$, continuous functions $M_1, M_2 : (\Re^+)^3 \to \Re^+$ and $\tilde{\gamma} \in N$ such that the following estimates hold for all $\xi \in \pi(t_0, x_0, u, d)$ and $t \geq \xi$:

$$\tilde{\varphi}(\gamma_1(\delta_1(t)V_2(t))) \leq \max \left\{ \begin{array}{l} S_2(M_2(t_0, \xi - t_0, \|x_0\|_X), t - \xi) \\ \tilde{\varphi}\left(\varphi^{-1}\left(\sup_{\xi \leq \tau \leq t} V_1(\tau)\right)\right), \sup_{t_0 \leq \tau \leq t} \tilde{\gamma}(\delta(\tau)\|u(\tau)\|_U) \end{array} \right\} \quad (3.37)$$

$$V_1(t) \leq \max \left\{ \begin{array}{l} S_1(M_1(t_0, \xi - t_0, \|x_0\|_X), t - \xi) \\ \sup_{\xi \leq \tau \leq t} \tilde{\varphi}(\gamma_1(\delta_1(\tau)V_2(\tau))), \sup_{t_0 \leq \tau \leq t} \tilde{\gamma}(\delta(\tau)\|u(\tau)\|_U) \end{array} \right\} \quad (3.38)$$

where $\delta \in K^+$ is defined by (3.9). Notice that if $\beta_1, \beta_2, c_1, c_2, \delta_2 \in K^+$ are bounded, then the functions $M_1, M_2 : (\Re^+)^3 \to \Re^+$ are independent of $t_0 \in \Re^+$ (but still depend on $\xi - t_0$). Moreover, the function $\tilde{\gamma} \in N$ in addition to (3.28) satisfies for all $s \geq 0$:

$$\tilde{\gamma}(s) \geq \max\left\{ \tilde{\kappa}(\kappa(\sigma_1(p_1^u(s), 0))), \tilde{\kappa}(\kappa(\sigma_1(p_1(\varphi(D(s))), 0))), \tilde{\kappa}(\varphi(\gamma_1^u(s))) \right\} \quad (3.39a)$$

$$\tilde{\gamma}(s) \geq \max\left\{ \tilde{\varphi}(\gamma_1(M\kappa(\varphi(\gamma_2^u(s))))), \tilde{\varphi}(\gamma_1(M\kappa(\kappa(\sigma_2(p_2^u(s), 0))))), \tilde{\varphi}(\gamma_1(M\kappa(\kappa(\sigma_2(p_2(\varphi(B(s))), 0))))) \right\} \quad (3.39b)$$

where $D(s) := \gamma_2^u(s) + \gamma_2\left(\varphi\left(\frac{1}{2}\varphi^2(\kappa(\varphi(\gamma_1^u(s)))) + \gamma_1(M\kappa(\varphi(\gamma_2^u(s))))\right)\right)$ and $B(s) := \kappa(\varphi(\gamma_1^u(s))) + \gamma_1(M\kappa(\varphi(\gamma_2^u(s))))$ are functions of class $N$.

We define:

$$a_1(h, T, R) := \sup\left\{ V_1(t_0 + h) - \sup_{t_0 \leq \tau \leq t_0 + h} \tilde{\gamma}(\delta(\tau)\|u(\tau)\|_U); \|x_0\|_X \leq R, t_0 \in [0, T], d \in M_D, u \in M_U \right\} \quad (3.40)$$

$$a_2(h, T, R) := \sup\left\{ \tilde{\varphi}(\gamma_1(\delta_1(t_0 + h)V_2(t_0 + h))) - \sup_{t_0 \leq \tau \leq t_0 + h} \tilde{\gamma}(\delta(\tau)\|u(\tau)\|_U); \|x_0\|_X \leq R, t_0 \in [0, T], d \in M_D, u \in M_U \right\} \quad (3.41)$$

$$l_1 := \limsup_{h \to +\infty} a_1(h, T, R), \quad l_2 := \limsup_{h \to +\infty} a_2(h, T, R) \quad (3.42)$$

By virtue of (3.24) and (3.27), the limits defined in (3.42) exist and are finite. Definition (3.42) implies that for every $\varepsilon > 0$, $T \geq 0$ and $R \geq 0$, there exists a $\tau := \tau(\varepsilon, T, R) \geq 0$, such that:

$$a_1(h, T, R) \leq l_1 + \varepsilon; \quad a_2(h, T, R) \leq l_2 + \varepsilon, \quad \forall h \geq \tau \quad (3.43)$$

By virtue of the Weak Semigroup Property for system $\Sigma$, there exists a constant $r > 0$, such that for each $(t_0, x_0, u, d) \in \Re^+ \times X \times M_U \times M_D$ we have $\pi(t_0, x_0, u, d) \cap [t_0 + \tau, t_0 + \tau + r] \neq \emptyset$. Let $\xi \in \pi(t_0, x_0, u, d) \cap [t_0 + \tau, t_0 + \tau + r]$. Estimates (3.37), (3.38) in conjunction with definitions (3.40), (3.41) and inequalities (3.43) give:

$$V_1(t) - \sup_{t_0 \leq \tau \leq t} \tilde{\gamma}(\delta(\tau)\|u(\tau)\|_U) \leq \max\left\{ S_1(M_1(t_0, \xi - t_0, \|x_0\|_X), t - \xi); l_2 + \varepsilon \right\} \quad (3.44)$$

$$\widetilde{\varphi}(\gamma_1(\delta_1(t)V_2(t))) - \sup_{t_0 \leq \tau \leq t} \widetilde{\gamma}\big(\delta(\tau)\|u(\tau)\|_{\mathcal{U}}\big) \leq \max\left\{ \begin{array}{l} S_2\big(M_2\big(t_0, \xi - t_0, \|x_0\|_X\big), t - \xi\big) \\ \widetilde{\varphi}\left(\varphi^{-1}\left(l_1 + \varepsilon + \sup_{t_0 \leq \tau \leq t} \widetilde{\gamma}\big(\delta(\tau)\|u(\tau)\|_{\mathcal{U}}\big)\right)\right) - \sup_{t_0 \leq \tau \leq t} \widetilde{\gamma}\big(\delta(\tau)\|u(\tau)\|_{\mathcal{U}}\big) \end{array} \right\}$$

(3.45)

Using the identity $\widetilde{\varphi}\big(\varphi^{-1}(s)\big) = s - \frac{1}{2}\rho\big(\varphi^{-1}(s)\big)$ and inequality (3.45) we get:

$$\widetilde{\varphi}(\gamma_1(\delta_1(t)V_2(t))) - \sup_{t_0 \leq \tau \leq t} \widetilde{\gamma}\big(\delta(\tau)\|u(\tau)\|_{\mathcal{U}}\big) \leq \max\left\{ S_2\big(M_2\big(t_0, \xi - t_0, \|x_0\|_X\big), t - \xi\big); l_1 + \varepsilon - \frac{1}{2}\rho\big(\varphi^{-1}(l_1)\big) \right\} \quad (3.46)$$

The properties of the $KL$ functions in conjunction with estimates (3.44), (3.46), the fact that $\xi \in [t_0 + \tau, t_0 + \tau + r]$ and definitions (3.40), (3.41), (3.42), give for all $\varepsilon > 0$:

$$l_1 \leq l_2 + \varepsilon \ ; \ l_2 \leq l_1 + \varepsilon - \frac{1}{2}\rho\big(\varphi^{-1}(l_1)\big)$$

From the first inequality we obtain $l_1 \leq l_2$. The second inequality implies $\rho\big(\varphi^{-1}(l_1)\big) \leq 2\varepsilon$ for all $\varepsilon > 0$, which directly gives $l_1 = l_2 = 0$.

Definitions (3.40), (3.41) and (3.42) imply that P3 of Lemma 2.16 holds for system $\Sigma$ with $V = V_1$ or $V = \widetilde{\varphi}\big(\gamma_1(\delta_1(t)V_2)\big)$, for appropriate $\widetilde{\gamma} \in N$ (which satisfies (3.28) and (3.39)) and $\delta \in K^+$ as defined by (3.9).

Notice that if $\beta_1, \beta_2, c_1, c_2, \delta_2 \in K^+$ are bounded then definitions (3.40), (3.41) are modified as follows:

$$a_1(h, R) := \sup\left\{ V_1(t_0 + h) - \sup_{t_0 \leq \tau \leq t_0 + h} \widetilde{\gamma}\big(\delta(\tau)\|u(\tau)\|_{\mathcal{U}}\big); \|x_0\|_X \leq R, t_0 \geq 0, d \in M_D, u \in M_U \right\} \quad (3.47)$$

$$a_2(h, R) :=$$
$$\sup\left\{ \varphi\big(\gamma_1(\delta_1(t_0 + h)V_2(t_0 + h))\big) - \sup_{t_0 \leq \tau \leq t_0 + h} \widetilde{\gamma}\big(\delta(\tau)\|u(\tau)\|_{\mathcal{U}}\big); \|x_0\|_X \leq R, t_0 \geq 0, d \in M_D, u \in M_U \right\} \quad (3.48)$$

Similar arguments as above show that property P3 of Lemma 2.17 holds for system $\Sigma$ with $V = V_1$ or $V = \widetilde{\varphi}\big(\gamma_1(\delta_1(t)V_2)\big)$, for appropriate $\widetilde{\gamma} \in N$ (which satisfies (3.28) and (3.39)) and $\delta \in K^+$ as defined by (3.9). The proof is complete. ◁

**Remark 3.6:** If the functions $\gamma_1^u, p_1^u, \gamma_2^u, p_2^u \in N$ are all identically zero then it follows that the gain function $\gamma \in N$ is identically zero. Indeed, the reader should notice that $\gamma \in N$ may be selected as $\gamma(s) := a\big(4\widetilde{\gamma}(s)\big) \in N$, where $a \in N$ is the function involved in hypothesis (H4) and $\widetilde{\gamma} \in N$ is the function that satisfies (3.28), (3.39a,b). Moreover, notice that for the input free case Theorem 3.1 and Corollary 3.4 imply (Uniform) Robust Global Asymptotic Output Stability (RGAOS) for the corresponding system. The following example shows the applicability of this particular remark to systems with impulses at fixed times.

**Example 3.7:** Consider the following system:

$$\dot{z}(t) = Az(t) + g(x(t)) \tag{3.49a}$$

$$\dot{x}(t) = f(x(t)), t \notin \pi$$

$$x(\tau_i) = h\left(\lim_{t \to \tau_i^-} x(t)\right) \tag{3.49b}$$

$$z(t) \in \Re^k, x(t) \in \Re^n$$

where $A \in \Re^{k \times k}$ is a Hurwitz matrix, $\pi = \{\tau_i\}_{i=0}^{\infty}$ is a partition of $\Re^+$ with diameter $r > 0$, $f : \Re^n \to \Re^n$, $g : \Re^n \to \Re^k$, $h : \Re^n \to \Re^n$ are continuous vector fields, $f(x)$ being locally Lipschitz with respect to $x \in \Re^n$, with $f(0) = 0$, $g(0) = 0$, $h(0) = 0$. Notice that subsystem (3.49a) is a system described by ODEs which satisfies hypotheses (A1-3) of Example 2.8. Hence, subsystem (3.49a) satisfies the BIC property and $0 \in \Re^k$ is a robust equilibrium point from the input $x$. Moreover, subsystem (3.49b) is a hybrid system with impulses at fixed times, which satisfies hypotheses (Q1-4) of Example 2.12. Hence, subsystem (3.49b) satisfies the BIC property and $0 \in \Re^n$ is a robust equilibrium point (from the zero input). We remark that since both subsystems (3.49a,b) satisfy the classical semigroup property and consequently the composite system (3.49) can be regarded as the feedback interconnection of subsystems (3.49a,b).

Since $A \in \Re^{k \times k}$ is Hurwitz, it follows that subsystem (3.49a) satisfies the UISS property from the input $x$. Moreover, if there exists a $C^1$ positive definite and radially unbounded function $V : \Re^n \to \Re^+$ and constants $c_1, c_2 \in \Re$ with $c_2 \neq 0$, $\mu, \lambda > 0$, such that:

$$\nabla V(x) f(x) \leq -c_1 V(x), \ \forall x \in \Re^n \tag{3.50a}$$

$$V(h(x)) \leq \exp(-c_2) V(x), \ \forall x \in \Re^n \tag{3.50b}$$

$$-c_2 \, card\big(\pi \cap [s, s+t)\big) \leq \mu + (c_1 - \lambda) t \ \forall s, t \in \Re^+ \tag{3.50c}$$

where $card(S)$ denotes the cardinal number of the set $S$, then Theorem 1 in [5] implies that $0 \in \Re^n$ is GAS for subsystem (3.49b). Taking into account Remark 3.2(b) and Remark 3.6, we conclude that $0 \in \Re^k \times \Re^n$ is uniformly GAS for the composite system (3.49) under the hypotheses stated above. ◁

## 4. Application to Partial State Sampled-Data Control

In this section we present applications of the Small-Gain results (Theorem 3.1 and Corollary 3.4) to partial state sampled-data control problems. It should be emphasized that sampled-data control systems cannot be handled with Small-Gain results that have appeared so far in the literature, since sampled-data control systems do not satisfy the classical semigroup property (see Example 2.11).

Consider the following control system described by ODEs:

$$\dot{z} = f(t, d, z, x, u)$$
$$z \in \Re^k, \ d \in D, u \in U, t \geq 0 \tag{4.1a}$$

$$\dot{x} = Ax + Bv + Bg(t, d, z, u)$$
$$x \in \Re^n, v \in \Re, d \in D, u \in U, t \geq 0 \tag{4.1b}$$

where $(A, B)$ is a controllable pair of matrices, $D \subset \Re^l$ is a compact set, $U \subseteq \Re^p$ is non-empty with $0 \in U$ and the mappings $f : \Re^+ \times D \times \Re^k \times \Re^n \times U \to \Re^k$, $g : \Re^+ \times D \times \Re^k \times U \to \Re$ are continuous, locally Lipschitz in $(z, x)$, uniformly in $d \in D$ with $f(t, d, 0, 0, 0) = 0$, $g(t, d, 0, 0) = 0$ for all $(t, d) \in \Re^+ \times D$. The problem we consider is the (W)ISS stabilization problem for (4.1) with sampled-data feedback applied with zero order hold and depending only on $x \in \Re^n$, i.e., we want to find a function $k : \Re^n \to \Re$ with $k(0) = 0$ and a constant $r > 0$ such that system (4.1a) with

$$\dot{x}(t) = Ax(t) + Bk(x(\tau_i)) + Bg(t, d(t), z(t), u(t)), \ t \in [\tau_i, \tau_{i+1})$$
$$\tau_{i+1} = \tau_i + \exp(-w(\tau_i)) r \ , \ w(t) \in \Re^+ \tag{4.2}$$

satisfies the WISS property from the inputs $(u, w)$. Notice that the input $w$ has been introduced in order to quantify the uncertainty in sampling times, i.e., we have to guarantee stability properties for the closed-loop system (4.1a)-(4.2) for all sampling schedules of diameter less than or equal to $r > 0$. To this purpose we make the following assumptions:

**(A1)** System (4.1a) satisfies the WISS property from the inputs $x$ and $u$. Specifically, there exist functions $\sigma \in KL$, $\beta, \delta_1^u \in K^+$, $\gamma_1, \gamma_1^u \in N$ such that for all $(t_0, z_0, d, x, u) \in \Re^+ \times \Re^k \times L_{loc}^\infty(\Re^+; D) \times L_{loc}^\infty(\Re^+; \Re^n) \times L_{loc}^\infty(\Re^+; U)$ the solution of (4.1) with initial condition $z(t_0) = z_0$ corresponding to inputs $(d, x, u) \in L_{loc}^\infty(\Re^+; D) \times L_{loc}^\infty(\Re^+; \Re^n) \times L_{loc}^\infty(\Re^+; U)$ satisfies the following estimate for all $t \geq t_0$:

$$|z(t)| \leq \sigma(\beta(t_0)|z_0|, t - t_0) + \sup_{t_0 \leq \tau \leq t} \gamma_1(|x(\tau)|) + \sup_{t_0 \leq \tau \leq t} \gamma_1^u(\delta_1^u(\tau)|u(\tau)|) \tag{4.3}$$

**(A2)** There exist functions $\delta_2^u \in K^+$, $\gamma_2, \gamma_2^u \in N$ such that the following inequality holds for all $(t, z, d, u) \in \Re^+ \times \Re^k \times D \times U$:

$$|g(t, d, z, u)| \leq \gamma_2(|z|) + \gamma_2^u(\delta_2^u(t)|u|) \tag{4.4}$$

**(A3)** There exist a function $\rho \in K_\infty$ and a constant $R \geq 1$ such that:

$$\gamma_1(R^{-1}\gamma_2(s) + \rho(R^{-1}\gamma_2(s))) + \rho(\gamma_1(R^{-1}\gamma_2(s) + \rho(R^{-1}\gamma_2(s)))) \leq s, \forall s \geq 0 \tag{4.5}$$

For example, hypothesis (A3) holds if $\gamma_i(s) = K_i s$, where $K_i \geq 0$ ($i = 1, 2$), i.e., if the gain functions are linear.

Next we show that the problem of WISS stabilization problem for (4.1) with sampled-data feedback applied with zero order hold and depending only on $x \in \Re^n$ is solvable under hypotheses (A1-3) by linear feedback. The proof of this result will be made by making use of Corollary 3.4.

Notice that since $(A, B)$ is a controllable pair of matrices, it follows that for every $\mu > 0$, $R \geq 1$ there exist a symmetric positive definite matrix $P \in \Re^{n \times n}$, a vector $k \in \Re^n$, constants $Q_1, Q_2 > 0$ such that the following inequalities hold for all $(x, u) \in \Re^n \times \Re$:

$$Q_1|x|^2 \leq x'Px \leq Q_2|x|^2$$
$$2x'P(A + Bk')x + 2x'PBu \leq -4\mu\, x'Px + \frac{Q_1}{4\mu R^2} u^2 \tag{4.6}$$

Next we show the following claim.

**Claim:** *For every $\mu > 0$, $R \geq 1$ there exists a vector $k \in \Re^n$, constants $M, r > 0$, such that for all $(t_0, x_0, d, z, u, w) \in \Re^+ \times \Re^n \times L_{loc}^\infty(\Re^+; D) \times L_{loc}^\infty(\Re^+; \Re^k) \times L_{loc}^\infty(\Re^+; U) \times L_{loc}^\infty(\Re^+; \Re^+)$ the solution of the hybrid system:*

$$\dot{x}(t) = Ax(t) + Bk'x(\tau_i) + Bg(t, d(t), z(t), u(t)), \; t \in [\tau_i, \tau_{i+1})$$
$$\tau_{i+1} = \tau_i + \exp(-w(\tau_i))\, r \quad, \quad w(t) \in \Re^+ \tag{4.7}$$

*with initial condition $x(t_0) = x_0$ corresponding to inputs $(d, z, u, w) \in L_{loc}^\infty(\Re^+; D) \times L_{loc}^\infty(\Re^+; \Re^k) \times L_{loc}^\infty(\Re^+; U) \times L_{loc}^\infty(\Re^+; \Re^+)$ satisfies the following estimate for all $t \geq t_0$:*

$$|x(t)| \leq M \exp(-\mu(t - t_0))|x_0| + \sup_{t_0 \leq \tau \leq t} R^{-1}\gamma_2(|z(\tau)|) + \sup_{t_0 \leq \tau \leq t} R^{-1}\gamma_2^u(\delta_2^u(\tau)|u(\tau)|) \tag{4.8}$$

*where $\delta_2^u \in K^+$, $\gamma_2, \gamma_2^u \in N$ are the functions involved in (4.4).*

**Proof of Claim:** Let arbitrary $\mu > 0$, $R > 1$. Since $(A, B)$ is a controllable pair of matrices, it follows that there exists a symmetric positive definite matrix $P \in \Re^{n \times n}$, a vector $k \in \Re^n$, constants $Q_1, Q_2 > 0$ such that inequalities (4.6) hold for all $(x, u) \in \Re^n \times \Re$. Let arbitrary $(t_0, x_0, d, z, u, w) \in \Re^+ \times \Re^n \times L^\infty_{loc}(\Re^+; D) \times L^\infty_{loc}(\Re^+; \Re^k) \times L^\infty_{loc}(\Re^+; U) \times L^\infty_{loc}(\Re^+; \Re^+)$ and consider the solution $x(t)$ of (4.7) with initial condition $x(t_0) = x_0$ corresponding to inputs $(d, z, u, w) \in L^\infty_{loc}(\Re^+; D) \times L^\infty_{loc}(\Re^+; \Re^k) \times L^\infty_{loc}(\Re^+; U) \times L^\infty_{loc}(\Re^+; \Re^+)$ (the solution exists for all $t \geq t_0$). Finally, consider the function $V(t) := x'(t) P x(t)$, which is absolutely continuous on $[t_0, +\infty)$. By virtue of (4.6) the derivative of $V(t)$ satisfies a.e. on the interval $[\tau_i, \tau_{i+1})$:

$$\dot{V}(t) \leq -4\mu V(t) + 2x'PBk'(x(\tau_i) - x(t)) + \frac{Q_1}{4\mu R^2} |g(t, d(t), z(t), u(t))|^2 \tag{4.9}$$

Let $r > 0$ constant that satisfies:

$$r \leq \frac{2\mu}{M|Bk'||A + Bk'|(2|A + Bk'| + |B|) + 2\mu|A| + 2\mu|A + Bk'|} \; ; \; r \leq \frac{1}{4\mu R^2 M|B||Bk'| + |A| + |A + Bk'|} \tag{4.10}$$

where $M := \frac{Q_2}{Q_1} \geq 1$.

It follows from (4.7) that $|x(t) - x(\tau_i)| \leq r|A| \sup_{\tau_i \leq s \leq t} |x(s) - x(\tau_i)| + r|A + Bk'||x(\tau_i)| + r|B| \sup_{\tau_i \leq s \leq t} |g(s, d(s), z(s), u(s))|$, which directly implies $|x(t) - x(\tau_i)| \leq \frac{r|A + Bk'|}{1 - r|A|} |x(\tau_i)| + \frac{r|B|}{1 - r|A|} \sup_{\tau_i \leq s \leq t} |g(s, d(s), z(s), u(s))|$, for all $t \in [\tau_i, \tau_{i+1})$. Moreover, the previous inequality in conjunction with the triangle inequality $|x(\tau_i)| \leq |x(t) - x(\tau_i)| + |x(t)|$ implies the estimate $|x(t) - x(\tau_i)| \leq \frac{r|A + Bk'|}{1 - r|A| - r|A + Bk'|} |x(t)| + \frac{r|B|}{1 - r|A| - r|A + Bk'|} \sup_{\tau_i \leq s \leq t} |g(s, d(s), z(s), u(s))|$, for all $t \in [\tau_i, \tau_{i+1})$. Using the previous inequality in conjunction with (4.9) and completing the squares, we obtain for almost all $t \in [\tau_i, \tau_{i+1})$:

$$\dot{V}(t) \leq -\left(4\mu - \frac{rM|Bk'||A + Bk'|(2|A + Bk'| + |B|)}{1 - r(|A| + |A + Bk'|)}\right) V(t) + \left(\frac{Q_1}{4\mu R^2} + \frac{r|B|Q_2|Bk'|}{1 - r(|A| + |A + Bk'|)}\right) \sup_{\tau_i \leq s \leq t} |g(s, d(s), z(s), u(s))|^2 \tag{4.11}$$

It follows from inequalities (4.10), (4.11) that the following estimate holds for the derivative of $V(t)$ a.e. on the interval $[\tau_i, \tau_{i+1})$:

$$\dot{V}(t) \leq -2\mu V(t) + \frac{Q_1}{2\mu R^2} \sup_{t_0 \leq s \leq t} |g(s, d(s), z(s), u(s))|^2 \tag{4.12}$$

Notice that since estimate (4.12) does not depend on the particular interval $[\tau_i, \tau_{i+1})$, we may conclude that estimate (4.12) holds a.e. for $t \geq t_0$. Estimate (4.12) implies directly that $V(t) \leq \exp(-2\mu(t - t_0)) V(t_0) + \frac{Q_1}{R^2} \sup_{t_0 \leq s \leq t} |g(s, d(s), z(s), u(s))|^2$ for all $t \geq t_0$. Finally, estimate (4.8) is an immediate consequence of the previous inequality, definitions $V(t) := x'(t) P x(t)$ and $M := \frac{Q_2}{Q_1} \geq 1$ as well as inequalities (4.4) and (4.6). The proof of the claim is complete. ◁

By making use of Corollary 3.4 and specifically Remark 3.5(b), we may conclude that the closed-loop system (4.1a) with (4.7) satisfies the WISS property from the inputs $u$ and $w$, when $R > 1$ is chosen to be greater than or equal to the constant involved in hypothesis (A3). Moreover the gain function for the input $w$ is identically zero. Furthermore, if the functions $\delta_1^u, \delta_2^u \in K^+$ are bounded, then the closed-loop system (4.1a) with (4.7) satisfies the ISS property from the inputs $u$ and $w$. Finally, if in addition $\beta \in K^+$ is bounded, then the closed-loop system (4.1a) with (4.7) satisfies the UISS property from the inputs $u$ and $w$.

**Example 4.1:** The following planar system described by ODEs:

$$\begin{aligned} \dot{z} &= -z^3 + z\,x \\ \dot{x} &= d(t)z^2 + u + v \\ (z,x)' &\in \Re^2, u, v \in \Re \end{aligned} \quad (4.13)$$

is studied in [9] where it is shown that if $d(t) \equiv d$ with $|d| < \frac{1}{2}$ then the feedback law $v(t) = -x(t)$, guarantees the UISS property for the closed-loop system from the input $u \in \Re$. The proof of this fact is made by using a slightly modified version of the Small-Gain Theorem presented in [10]. Here we study the possibility of robustly globally stabilizing the origin for system (4.13), using the following feedback law with zero order hold and positive sampling rate:

$$\begin{aligned} v(t) &= -x(\tau_i), t \in [\tau_i, \tau_{i+1}) \\ \tau_{i+1} &= \tau_i + \exp(-w(\tau_i))\,r, w(t) \in \Re^+ \end{aligned} \quad (4.14)$$

for time-varying disturbances $d(t) \in D := [-\delta, \delta]$ with $\delta \in (0,1)$.

First notice that system (4.13) is a system of the form (4.1) with $f(t,d,z,x,u) = -z^3 + z\,x$, $g(t,d,z,u) = d\,z^2 + u$, $B = [1]$, $A = [0]$. Moreover, hypothesis (A2) holds with $\gamma_2(s) = \delta s^2$, $\gamma_2^u(s) := s$ and $\delta_2^u(t) \equiv 1$.

Working exactly as in [9] it may be shown that for every $\varepsilon \in (0,1)$ the subsystem $\Sigma_1$:

$$\dot{z} = -z^3 + z\,x \quad (4.15)$$

satisfies the UISS property from the input $x$ with gain function $\gamma_1(s) := \sqrt{\dfrac{s}{1-\varepsilon}}$. Thus the subsystem $\Sigma_1$ satisfies hypothesis (A1) with $\gamma_1(s) := \sqrt{\dfrac{s}{1-\varepsilon}}$, $\beta(t) \equiv 1$, $\gamma_1^u(s) \equiv 0$ and appropriate $\sigma \in KL$.

For every $\delta \in (0,1)$ there exist $\varepsilon \in (0,1)$ and $L > 0$ such that:

$$(1+L)\sqrt{\dfrac{(1+L)\delta}{1-\varepsilon}} \leq 1 \quad (4.16)$$

Selecting $\rho(s) := Ls$ with $L > 0$, we conclude from (4.16) that (4.5) holds with $R = 1$. Finally, since (4.6) holds with $Q_1 = Q_2 = 1$, $P = [1]$, $\mu = \dfrac{1}{4}$ and $k = -1$, we conclude that (4.13) with (4.14) satisfies the UISS property from the inputs $u$ and $w$. Moreover the gain function for the input $w$ is identically zero. The maximum allowable sampling period ($r$) may be determined by inequalities (4.10), which give $r = 1/7$. ◁

# 5. Conclusions

A Small-Gain Theorem, which can be applied to a wide class of systems that includes systems that satisfy the weak semigroup property, is presented in the present work. The result generalizes all existing results in the literature and exploits notions of weighted, uniform and non-uniform Input-to-Output Stability (IOS) property. Moreover, the Small-Gain Theorem of the present work is a method for establishing qualitative properties expressed in a very general framework unifying works from various fields as well as different stability notions. The results presented in the paper can be extended without much difficulty to the case of local stability notions.

Applications to partial state feedback stabilization problems with sampled-data feedback applied with zero order hold and positive sampling rate, are also presented. It should be emphasized that sampled-data control systems cannot be handled with Small-Gain results that have appeared so far in the literature, since sampled-data control systems do not satisfy the classical semigroup property. The results are illustrated by examples, which show the usefulness of the main result for the stability analysis of interconnected systems.

**Acknowledgments:** This work is supported partly by the U.S. NSF under grants ECS-0093176, OISE-0408925 and DMS-0504462.

# Appendix

**Proof of Lemma 2.13:** Lemma 3.5 in [20] guarantees that the control system $\Sigma := (X, Y, M_U, M_D, \phi, \pi, H)$ which has the BIC property is RFC from the input $u \in M(U)$ if and only if there exist functions $q \in K^+$, $a \in K_\infty$ and a constant $R \geq 0$ such that the following estimate holds for all $(t_0, x_0, d, u) \in \Re^+ \times X \times M_D \times M_U$ and $t \geq t_0$:

$$\|\phi(t, t_0, x_0, u, d)\|_X \leq q(t) a\left( R + \|x_0\|_X + \sup_{t_0 \leq \tau \leq t} \|u(\tau)\|_U \right) \tag{A1}$$

It should be emphasized that all results in [20] were proved under the assumption of the classical semigroup property for the control system. However, the proof of Lemma 3.5 does not depend on the semigroup property and consequently may be repeated as it stands for a system, which satisfies the weak semigroup property. Let $\beta \in K^+$ arbitrary. Using (A1) we obtain:

$$\begin{aligned}
\beta(t)\|\phi(t, t_0, x_0, u, d)\|_X &\leq \frac{1}{2} q^2(t)\beta^2(t) + \frac{1}{2}\max\left\{ a^2(3R); a^2(3\|x_0\|_X); \sup_{t_0 \leq \tau \leq t} a^2(3\|u(\tau)\|_U) \right\} \\
&\leq \max\left\{ q^2(t)\beta^2(t); a^2(3R); a^2(3\|x_0\|_X); \sup_{t_0 \leq \tau \leq t} a^2(3\|u(\tau)\|_U) \right\} \\
&\leq \max\left\{ \gamma(t); a^2(3\|x_0\|_X); \sup_{t_0 \leq \tau \leq t} a^2(3\|u(\tau)\|_U) \right\}
\end{aligned} \tag{A2}$$

where $\gamma(t) = q^2(t)\beta^2(t) + a^2(3R)$. Define:

$$a(T, s) := \max\{\gamma(t_0 + h) - \gamma(t_0) : h \in [0, s], t_0 \in [0, T]\} \tag{A3}$$

Clearly, definition (A3) implies that for each fixed $s \geq 0$ $a(\cdot, s)$ is non-decreasing and for each fixed $T \geq 0$ $a(T, \cdot)$ is non-decreasing. Furthermore, continuity of $\gamma$ guarantees that for every $T \geq 0$ $\lim_{s \to 0^+} a(T, s) = a(T, 0) = 0$. It turns out from Lemma 2.3 in [17], that there exist functions $\zeta \in K_\infty$ and $\kappa \in K^+$ such that

$$a(T, s) \leq \zeta(\kappa(T)s), \quad \forall (T, s) \in (\Re^+)^2 \tag{A4}$$

Combining definition (A3) with inequality (A4), we conclude that for all $t_0 \geq 0$ and $t \geq t_0$, it holds that:

$$\begin{aligned}
\gamma(t) &\leq \gamma(t_0) + \zeta(\kappa(t_0)(t - t_0)) \leq \gamma(t_0) + \zeta\left(\frac{1}{2}\kappa^2(t_0) + \frac{1}{2}(t - t_0)^2\right) \\
&\leq \gamma(t_0) + \zeta(\kappa^2(t_0)) + \zeta((t - t_0)^2) \leq \max\{2\gamma(t_0) + 2\zeta(\kappa^2(t_0)); 2\zeta((t - t_0)^2)\}
\end{aligned}$$

The above inequality in conjunction with (A2) implies that (2.9) holds for all $(t_0, x_0, d, u) \in \Re^+ \times X \times M_D \times M_U$ and $t \geq t_0$ with $\mu(t) := 2\zeta(t^2)$, $c(t) := 2\gamma(t) + 2\zeta(\kappa^2(t))$, $a(s) := a^2(3s)$ and $p(s) := a^2(3s)$. The proof is complete. ◁

**Proof of Lemma 2.16:** As in the proof of Proposition 2.2 in [17], let $T, h \geq 0$, $s \geq 0$ and define:

$$a(T, s) := \sup\left\{ V(t, \phi(t, t_0, x_0, u, d), u(t)) - \sup_{t_0 \leq \tau \leq t} \gamma(\delta(\tau)\|u(\tau)\|_U); \|x_0\|_X \leq s, t \geq t_0 \in [0, T], d \in M_D, u \in M_U \right\} \tag{A5}$$

$$M(h, T, s) := \sup\left\{ V(t_0 + h, \phi(t_0 + h, t_0, x_0, u, d), u(t)) - \sup_{t_0 \leq \tau \leq t_0 + h} \gamma(\delta(\tau)\|u(\tau)\|_U); \|x_0\|_X \leq s, t_0 \in [0, T], d \in M_D, u \in M_U \right\} \tag{A6}$$

First notice that by virtue of property P1 it holds that $a(T,s) < +\infty$ for all $T \geq 0$, $s \geq 0$. Moreover, notice that since $0 \in X$ is a robust equilibrium point from the input $u \in M_U$ and $V(t,0,0) = 0$ for all $t \geq 0$, we have $a(T,s) \geq 0$ for all $T \geq 0$, $s \geq 0$. Furthermore, notice that $M$ is well-defined, since by definitions (A5), (A6) the following inequality is satisfied for all $T, h \geq 0$ and $s \geq 0$:

$$0 \leq M(h,T,s) \leq a(T,s) \tag{A7}$$

Clearly, definition (A5) implies that for each fixed $s \geq 0$ $a(\cdot, s)$ is non-decreasing and for each fixed $T \geq 0$ $a(T, \cdot)$ is non-decreasing. Furthermore, property P2 asserts that for every $T \geq 0$ $\lim_{s \to 0^+} a(T,s) = 0$. Hence, the inequality $a(T,0) \geq 0$ for all $T \geq 0$, in conjunction with $\lim_{s \to 0^+} a(T,s) = 0$ and the fact that $a(T, \cdot)$ is non-decreasing implies $a(\cdot, 0) = 0$. It turns out from Lemma 2.3 in [17], that there exist functions $\zeta \in K_\infty$ and $q \in K^+$ such that

$$a(T,s) \leq \zeta(q(T)s), \quad \forall (T,s) \in (\Re^+)^2 \tag{A8}$$

Without loss of generality we may assume that $q \in K^+$ is non-decreasing. Moreover, property P3 guarantees that for every $\varepsilon > 0$, $T \geq 0$ and $R \geq 0$, there exists a $\tau = \tau(\varepsilon, T, R) \geq 0$, such that

$$M(h,T,s) \leq \varepsilon \text{ for all } h \geq \tau(\varepsilon,T,R) \text{ and } 0 \leq s \leq R \tag{A9}$$

Let

$$g(s) := \sqrt{s} + s^2 \tag{A10}$$

and let $p$ be a non-decreasing function of class $K^+$ with $p(0) = 1$ and

$$\lim_{t \to +\infty} p(t) = +\infty \tag{A11}$$

Define

$$\mu(h) := \sup\left\{ \frac{M(h,T,s)}{p(T)g(\zeta(q(T)s))}; \quad T \geq 0, s > 0 \right\} \tag{A12}$$

Obviously, by virtue of (A7), (A8) and (A10) the function $\mu : \Re^+ \to \Re^+$ is well defined and satisfies $\mu(\cdot) \leq 1$. We show that $\lim_{h \to +\infty} \mu(h) = 0$, equivalently, we establish that for any given $\varepsilon > 0$, there exists a $\delta = \delta(\varepsilon) \geq 0$ such that

$$\mu(h) \leq \varepsilon, \text{ for } h \geq \delta(\varepsilon) \tag{A13}$$

Notice first, that for any given $\varepsilon > 0$ there exist constants $a := a(\varepsilon)$ and $b := b(\varepsilon)$ with $0 < a < b$ such that

$$x \notin (a,b) \Rightarrow \frac{x}{\sqrt{x} + x^2} \leq \varepsilon \tag{A14}$$

We next recall (A11), which asserts that, for the above $\varepsilon$ for which (A14) holds, there exists a $c := c(\varepsilon) \geq 0$ such that $p(T) \geq \frac{1}{\varepsilon}$ for all $T \geq c$. This by virtue of (A7), (A8), (A10) and (A14) yields:

$$\frac{M(h,T,s)}{p(T)g(\zeta(q(T)s))} \leq \varepsilon, \quad \forall h \geq 0, \text{ when } T \geq c, \text{ or } \zeta(q(T)s) \notin (a,b) \tag{A15}$$

Hence, in order to establish (A13), it remains to consider the case:

$$a \leq \zeta(q(T)s) \leq b \text{ and } 0 \leq T \leq c \tag{A16}$$

Since, for each fixed $(h,s) \in (\Re^+)^2$ the mappings $M(h,\cdot,s)$, $M(h,s,\cdot)$, $q(\cdot)$ and $p(\cdot)$ are non-decreasing, we have that

$$\frac{M(h,T,s)}{p(T)g(\zeta(q(T)s))} \leq \frac{M\left(h,c,\frac{\zeta^{-1}(b)}{q(0)}\right)}{g(a)} \tag{A17}$$

provided that (A16) holds. By using (A9) and (A17) with

$$\varepsilon := \varepsilon\, g(a),\, T := c,\, R := \frac{\zeta^{-1}(b)}{q(0)}$$

it follows

$$M\left(h,c,\frac{\zeta^{-1}(b)}{q(0)}\right) \leq \varepsilon\, g(a),\text{ for } h \geq \delta(\varepsilon) := \tau\left(\varepsilon\, g(a), c, \frac{\zeta^{-1}(b)}{q(0)}\right) \tag{A18}$$

By taking into account (A15), (A16), (A17), (A18) and definition (A12) of $\mu(\cdot)$ it follows that (A13) holds with $\delta = \delta(\varepsilon)$ as selected in (A18). Since $\varepsilon > 0$ was arbitrary we conclude that $\lim_{h \to +\infty} \mu(h) = 0$. Consequently, there exists a continuous strictly decreasing function $\bar{\mu} : \Re^+ \to (0,+\infty)$ such that $\bar{\mu}(h) \geq \mu(h)$ for all $h \geq 0$ and $\lim_{h \to +\infty} \bar{\mu}(h) = 0$. Thus, by recalling definition (A12) we obtain

$$M(h,T,s) \leq \bar{\mu}(h)\theta(T,s),\; \forall (T,s) \in (\Re^+)^2,\; \forall h \geq 0 \tag{A19}$$

where $\theta(T,s) := p(T)g(\zeta(q(T)s))$. Clearly, $\theta$ satisfies all hypotheses of Lemma 2.3 in [17] and therefore there exist $\zeta_2 \in K_\infty$ and $\beta \in K^+$ such that

$$\theta(T,s) \leq \zeta_2(\beta(T)s),\; \forall (T,s) \in (\Re^+)^2 \tag{A20}$$

Thus definition (A6) implies that the following estimate holds for all $u \in M_U$, $(t_0, x_0, d) \in \Re^+ \times X \times M_D$ and $t \geq t_0$:

$$V(t,\phi(t,t_0,x_0,u,d),u(t)) \leq \bar{\mu}(t-t_0)\zeta_2\left(\beta(t_0)\|x_0\|_X\right) + \sup_{t_0 \leq \tau \leq t} \gamma\left(\delta(\tau)\|u(\tau)\|_U\right) \tag{A21}$$

Estimate (A21) implies (2.12) with $\sigma(s,t) := \bar{\mu}(t)\zeta_2(s)$. ◁

**Proof of Lemma 2.17:** As in the proof of Lemma 2.16, let $h \geq 0$, $s \geq 0$ and define:

$$a(s) := \sup\left\{V(t,\phi(t,t_0,x_0,u,d),u(t)) - \sup_{t_0 \leq \tau \leq t}\gamma\left(\delta(\tau)\|u(\tau)\|_U\right); \|x_0\|_X \leq s, t \geq t_0 \geq 0, d \in M_D, u \in M_U\right\} \tag{A22}$$

$$M(h,s) := \sup\left\{V(t_0+h,\phi(t_0+h,t_0,x_0,u,d),u(t)) - \sup_{t_0 \leq \tau \leq t_0+h}\gamma\left(\delta(\tau)\|u(\tau)\|_U\right); \|x_0\|_X \leq s, t_0 \geq 0, d \in M_D, u \in M_U\right\} \tag{A23}$$

First notice that by virtue of property P1 it holds that $a(s) < +\infty$ for all $s \geq 0$. Moreover, notice that since $0 \in X$ is a robust equilibrium point from the input $u \in M_U$ and $V(t,0,0) = 0$ for all $t \geq 0$, we have $a(s) \geq 0$ for all $s \geq 0$. Furthermore, notice that $M$ is well-defined, since by definitions (A22), (A23) the following inequality is satisfied for all $h \geq 0$ and $s \geq 0$:

$$0 \leq M(h,s) \leq a(s) \tag{A24}$$

Clearly, definition (A22) implies that $a(\cdot)$ is non-decreasing. Furthermore, property P2 asserts that $\lim_{s \to 0^+} a(s) = 0$. Hence, the inequality $a(0) \geq 0$, in conjunction with $\lim_{s \to 0^+} a(s) = 0$ and the fact that $a(\cdot)$ is non-decreasing implies $a(0) = 0$. It turns out that $a$ can be bounded from above by the $K_\infty$ function $\tilde{a}$ defined by $\tilde{a}(s) := s + \frac{1}{s} \int_s^{2s} a(w) dw$ for $s > 0$ and $\tilde{a}(0) = 0$. Define

$$\mu(h) := \sup\left\{ \frac{M(h,s)}{g(\tilde{a}(s))} ; s > 0 \right\} \tag{A25}$$

where $g$ is defined by (A10). Working exactly as in the proof of Lemma 2.16 we can show that the function $\mu : \Re^+ \to \Re^+$ is well defined and satisfies $\mu(\cdot) \leq 1$, $\lim_{h \to +\infty} \mu(h) = 0$. Consequently, there exists a continuous strictly decreasing function $\bar{\mu} : \Re^+ \to (0, +\infty)$ such that $\bar{\mu}(h) \geq \mu(h)$ for all $h \geq 0$ and $\lim_{h \to +\infty} \bar{\mu}(h) = 0$. Thus, by recalling definition (A25) we obtain

$$M(h,s) \leq \bar{\mu}(h) g(\tilde{a}(s)), \quad \forall h, s \geq 0 \tag{A26}$$

Hence definition (A23) implies that the following estimate holds for all $u \in M_U$, $(t_0, x_0, d) \in \Re^+ \times X \times M_D$ and $t \geq t_0$:

$$V(t, \phi(t, t_0, x_0, u, d), u(t)) \leq \bar{\mu}(t - t_0) g\left(\tilde{a}(\|x_0\|_X)\right) + \sup_{t_0 \leq \tau \leq t} \gamma\left(\delta(\tau) \|u(\tau)\|_U\right) \tag{A27}$$

Estimate (A27) implies (2.12) with $\beta(t) \equiv 1$ and $\sigma(s, t) := \bar{\mu}(t) g(\tilde{a}(s))$.  ◁

**Proof of Lemmas 2.18-2.19:** The proof is based on the following observation: if $\Sigma := (X, Y, M_U, M_D, \phi, \pi, H)$ is $T$-*periodic* then for all $(t_0, x_0, u, d) \in \Re^+ \times X \times M_U \times M_D$ it holds that $\phi(t, t_0, x_0, u, d) = \phi(t - kT, t_0 - kT, x_0, P_{kT}u, P_{kT}d)$ and $H(t, \phi(t, t_0, x_0, u, d), u(t)) = H(t - kT, \phi(t - kT, t_0 - kT, x_0, P_{kT}u, P_{kT}d), (P_{kT}u)(t - kT))$, where $k := [t_0 / T]$ denotes the integer part of $t_0 / T$ and the inputs $P_{kT}u \in M_U$, $P_{kT}d \in M_D$ are defined in Definition 2.2.

Since $\Sigma := (X, Y, M_U, M_D, \phi, \pi, H)$ satisfies the WIOS property from the input $u \in M_U$, there exist functions $\sigma \in KL$, $\beta, \delta \in K^+$, $\gamma \in N$ such that (2.10) holds for all $(t_0, x_0, u, d) \in \Re^+ \times X \times M_U \times M_D$ and $t \geq t_0$. Consequently, it follows that the following estimate holds for all $(t_0, x_0, u, d) \in \Re^+ \times X \times M_U \times M_D$ and $t \geq t_0$:

$$\|H(t, \phi(t, t_0, x_0, u, d), u(t))\|_Y \leq \sigma\left(\beta(t_0 - kT)\|x_0\|_X, t - t_0\right) + \sup_{\tau \in [t_0 - kT, t - kT]} \gamma\left(\delta(\tau) \|(P_{kT}u)(\tau)\|_U\right)$$

Setting $\tau = s - kT$ and since $0 \leq t_0 - \left[\frac{t_0}{T}\right]T < T$, for all $t_0 \geq 0$, we obtain

$$\|H(t, \phi(t, t_0, x_0, u, d), u(t))\|_Y \leq \tilde{\sigma}\left(\|x_0\|_X, t - t_0\right) + \sup_{s \in [t_0, t]} \gamma\left(\delta(s - kT) \|(P_{kT}u)(s - kT)\|_U\right) \tag{A28}$$

where $\tilde{\sigma}(s, t) := \sigma(rs, t)$ and $r := \max\{\beta(t); 0 \leq t \leq T\}$. Estimate (A28) and the identity $(P_{kT}u)(s - kT) = u(s)$ for all $s \geq 0$, implies that the following estimate holds for all $(t_0, x_0, u, d) \in \Re^+ \times X \times M_U \times M_D$ and $t \geq t_0$:

$$\|H(t, \phi(t, t_0, x_0, u, d), u(t))\|_Y \leq \tilde{\sigma}\left(\|x_0\|_X, t - t_0\right) + \sup_{s \in [t_0, t]} \gamma\left(\tilde{\delta}(s) \|u(s)\|_U\right) \tag{A29}$$

where $\tilde{\delta}(t) := \max\{\delta(s); s \in [0, t]\}$.

In case that $\Sigma := (X, Y, M_U, M_D, \phi, \pi, H)$ satisfies the IOS property from the input $u \in M_U$, then all arguments above may be repeated with $\delta(t) \equiv 1$. Thus we conclude that (A29) holds for all $(t_0, x_0, u, d) \in \Re^+ \times X \times M_U \times M_D$ and $t \geq t_0$ with $\tilde{\delta}(t) \equiv 1$. The proof is complete.  ◁